\documentclass{article}
\usepackage{latexsym}
\usepackage{amssymb,amsbsy,amsmath,amsfonts,amssymb,amscd}
\usepackage{graphicx,color}
\usepackage{float,url}
\usepackage{cancel}
\usepackage[colorlinks=true, citecolor =violet, linkcolor = violet, urlcolor  = violet]{hyperref}
\usepackage{xcolor}
\usepackage[notref,notcite]{showkeys}
\usepackage{tikz}
\usepackage{bbm}
\usepackage[normalem]{ulem}
\definecolor{relbleu}{RGB}{0,130,205}
\definecolor{relrose}{RGB}{214,0,120}
\definecolor{relorange}{RGB}{224,116,0}
\providecolor{bluegreen}{RGB}{0,139,139}

\definecolor{relvert}{RGB}{0,140,70}

\usepackage{amsthm}
\newtheorem{theorem}{Theorem}[section]

\newtheorem{lemma}[theorem]{Lemma}

\newtheorem{prop}[theorem]{Proposition}

\newtheorem{remark}[theorem]{Remark}

\newcommand{\Cc}{\mathcal{C}}

\newcommand{\N}{\mathbb{N}}

\newcommand{\R}{\mathbb{R}}

\newcommand{\eps}{\varepsilon}
\newcommand{\fa}{\forall}

\newcommand{\commentout}[1]{}

\title{On the Integrate-and-Fire model with superlinear drift
}

\date{\today}

\author{Beno\^\i t Perthame%
\thanks{Sorbonne Universit\'e, CNRS, Universit\'e de Paris Cit\'e, Inria,
Laboratoire Jacques-Louis Lions (LJLL), F-75005 Paris, France.
\newline\texttt{benoit.perthame@sorbonne-universite.fr}}
\and
Cl\'ement Rieutord%
\:\thanks{Corresponding author, Sorbonne Universit\'e, CNRS, Universit\'e de Paris Cit\'e,
Laboratoire Jacques-Louis Lions (LJLL), F-75005 Paris, France.
\newline\texttt{clement.rieutord@sorbonne-universite.fr}}
\and
Delphine Salort%
\thanks{Sorbonne Universit\'e, CNRS, Universit\'e de Paris Cit\'e,
Laboratoire Jacques-Louis Lions (LJLL), F-75005 Paris, France.
\newline\texttt{dsalort@gmail.com}}
}

\begin{document}

\maketitle

\begin{abstract}
The Nonlinear Integrate-and-Fire equation is a fundamental model for neural assemblies, describing the dynamics of the density of neurons with a given potential. Its standard formulation uses a given firing potential, and is known to undergo finite-time blow-up for excitatory connections. 

To fit experimental observations, it has been advocated to use superlinear drifts. We show how, in this situation, it makes sense to set the equation on the whole line. The main question we address is whether in these circumstances the solution still exhibits finite-time blow-up. 

We study the construction of solutions by a fixed point and, using the framework of \cite{carrilloQualitativePropertiesSolutions2015}, we develop novel a priori estimates which use specifically the superlinearity of the drift at $+\infty$. In particular, we obtain Cesàro-mean bounds on the activity, as well as exponential convergence towards the unique steady state in the weakly connected regime. Several critical steps rely on controls depending solely on  the $L^1$ norm of the solution or coefficients. 
\end{abstract}

\noindent{\makebox[1in]\hrulefill}\newline


\tableofcontents

\section*{Introduction}

Assemblies of Integrate-and-Fire neurons are best fit to experimental data when the potential drift is superlinear, see~\cite{bretteAdaptiveExponentialIntegrateandfire2005,brunelFiringRateNoisy2003}. This property only becomes mathematically meaningful when the membrane potential is assumed to evolve on the whole real line. Indeed, when the membrane potential  is upper bounded, the firing mechanism is naturally associated with reaching a finite threshold $V_F$ and superlinearity does not make sense. In contrast, on the full line, a superlinear drift provides an intrinsic firing mechanism by driving the potential to $+\infty$ in finite time. This dynamic can therefore be interpreted as the analogue of threshold crossing for an unbounded potential, as proposed by~\cite{perthameFokkerPlanckEquationSuperlinear2026} in the linear setting describing independent neurons. Here we investigate its nonlinear extension, corresponding to homogeneously connected neurons, which we call the Superlinear Integrate-and-Fire model (SLIF model).
\begin{equation}\label{eq:FPSL Non-lin}
    \left\{\begin{aligned}
        & \frac{\partial u}{\partial t}(t,x) + \frac{\partial}{\partial x}\big((h(x)+bN(t))u(t,x)\big) - \frac{\partial^2 u}{\partial x^2}(t,x) = N(t)\delta_{0}(x), \quad x \in \mathbb{R},\ t \geq 0,
        \\
        & N(t) := \lim_{x \rightarrow \infty} h(x)u(t,x),
        \\
        & u(0,x) = u^0(x).
    \end{aligned}\right.
\end{equation}
Here, $x$ denotes the membrane potential of neurons, $u(t,x)$ is the probability density of the membrane potential of a typical neuron in the network, $N(t)$ its firing rate, $h$ represents a superlinear drift, and $b\in\mathbb{R}$ denotes the mean connectivity of the network. The value $b>0$ corresponds to an excitatory network, while $b<0$ describes inhibitory interactions.

In comparison, in the classical setting, the firing mechanism is described by prescribing a finite firing threshold $V_F>0$ (with $0$ chosen as the reset potential). The corresponding nonlinear noisy leaky integrate-and-fire (NNLIF) equation reads
\begin{equation}\label{eq:NNLIF}
    \left\{\begin{aligned}
        & \frac{\partial u}{\partial t}(t,x) + \frac{\partial}{\partial x}\big((h(x)+bN(t))u(t,x)\big) - \frac{\partial^2 u}{\partial x^2}(t,x) = N(t)\delta_{0}(x), \quad x \leq V_F,\ t \geq 0,
        \\
        & u(t,V_F)=0, \qquad N(t):=-\frac{\partial u}{\partial x}(t,V_F).
    \end{aligned}\right.
\end{equation}
The above NNLIF model has been extensively studied and its different properties are well established (see the survey \cite{carrilloNonlinearPartialDifferential2025}). In particular, when $b>0$, solutions of the NNLIF equation may exhibit finite-time blow-up.  In the present work, we investigate whether this blow-up phenomenon persists when the finite threshold mechanism is replaced by the superlinear drift mechanism of Eq.~\eqref{eq:FPSL Non-lin}.
%, where neurons potential can reach $+\infty$ in finite time.

The possible blow-up mechanisms driving the dynamics of the two models can be compared. For the NNLIF equation, finite-time blow-up is characterized numerically by
\[
\limsup_{t\rightarrow T}N(t)=+\infty ,
\]
and theoretically by the fact that no solution of \eqref{eq:NNLIF} can be continued beyond $T$.
However, extensions after blow-up can be constructed using a dilated time scale, see~\cite{douNoisyIntegrateandfireEquation2024}, or through the probabilistic notion of physical solution, see~\cite{delarueParticleSystemsSingular2015, delarueGlobalSolvabilityNetworked2015, hamblyMcKeanVlasovEquationPositive2019}. The intuitive mechanism behind this phenomenon is the following. In the classical NNLIF model, the density $u(t)$ is supported on the half-line $(-\infty,V_F]$, and neurons fire when they reach the threshold $V_F$. When $b>0$, the instantaneous coupling through the firing rate modifies the drift by adding the term $bN(t)$, which pushes the population towards the firing boundary. As neurons accumulate near $V_F$, diffusion generates an increasing outgoing flux $N(t)$; this flux instantaneously feeds back into the drift, further accelerating the transport towards the threshold. The interplay between diffusion, threshold crossing, and instantaneous excitatory coupling creates a positive feedback loop that can lead to a divergence of the firing rate. In particular, for sufficiently strong excitatory interactions, the NNLIF model fails to admit stationary states, and this regime is associated with systematically finite-time blow-up.

A way to extend solutions while taking this phenomenon into account was first proposed in \cite{zbMATH08050032}, and later extended to a constant diffusion coefficient in \cite{douNoisyIntegrateandfireEquation2024}. The idea is to rescale time according to the firing rate, through the change of variable
\[
d\tau = N(t)\,dt
\]
(or $d\tau=(N(t)+c)\,dt$, $c>0$, as in \cite{zbMATH08050032}). In \cite{douNoisyIntegrateandfireEquation2024}, the NNLIF equation in these coordinates is obtained as the limit $\varepsilon\to0$ of a regularized version, in which the firing at the threshold $V_F$ is replaced by a random discharge with rate $\frac{1}{\varepsilon}$ on $\{x\ge V_F\}$. This method distinguishes the time regions where blow-up occurs from the periods where the solution is classical, while preserving a dynamics on both types of periods.
Let us mention \cite{liuInvestigatingIntegrateFire2021} as well, in which the NNLIF equation (in its linear regime) is rigorously obtained as a limit of a process of random discharge, in a similar spirit to \cite{douNoisyIntegrateandfireEquation2024}.

For the SLIF model, the situation is different. Here, the firing mechanism is no longer associated with reaching a prescribed finite threshold, but with the finite-time escape of probability mass towards $+\infty$ induced by the superlinear drift. The additional drift generated by the network activity, $bN(t)$, modifies this escape mechanism through a spatially homogeneous transport contribution, while the superlinear drift governs the behavior of the dynamics for large potentials. In contrast with the NNLIF model, the SLIF equation always admits stationary states. The main question addressed in this paper is therefore whether the feedback induced by the firing rate is strong enough to destabilize this stationary regime and still produce finite-time blow-up. 

We answer this question in the negative: the solution is global for every $b\in \R$.

\paragraph{Modeling of the neuron network.} 
%\label{sec:Model Neuro}
The membrane potential of a neuron, denoted $V_t$, is modulated by ion channels that block or allow the flow of ions between the inside and the outside of the membrane. It is therefore natural to view a neuron as an electrical circuit. However, when the potential approaches a threshold, its dynamics become strongly nonlinear and evolve on a fast time scale: the analogy with a simple electrical circuit no longer holds.

Mathematically speaking, to account for this dynamics specific to the neuron we let the membrane potential, which may now take arbitrarily large values, follow a stochastic differential equation whose drift $h$ has two prescribed asymptotic behaviors. 
At $-\infty$, in agreement with the NNLIF model, $h$ is linear,
\[
h(x)\approx -x,
\]
while at $+\infty$ it grows superlinearly (for instance $h(x) \sim_{x\to\infty} x^\alpha$ with $\alpha>1$, or $h(x) \sim_{x\to\infty} e^x$).

With such a drift, the trajectories $x(t)$ solving
\begin{equation}
    \left\{\begin{aligned}
        &\frac{d x(t)}{dt} = h(x(t)),\\
        & x(0)=x_0\gg1,
    \end{aligned}
    \right.
\end{equation}
reach $+\infty$ in finite time. Once this boundary is reached, the potential is instantaneously reset to $0$. The membrane potential then follows the dynamics
\begin{equation}\label{eq:EDS pot}
    \left\{\begin{aligned}
        & dV_t = h(V_t)\,dt+dI(t),        \\
        & \text{when} \ V_{t^-} = +\infty, \quad V_{t^+} = 0.
      \end{aligned}\right.
\end{equation}
The term $I(t)$ represents the incoming synaptic current. It describes the coupling of a neuron with the whole network, in the limit of a large number of neurons and under a mean-field interaction; its modeling relies on several assumptions, detailed in \cite{brunelDynamicsSparselyConnected, mattiaPopulationDynamicsInteracting2002, ricciardiOrnsteinUhlenbeckProcessModel1979}. Its rigorous derivation from a system of $K$ interacting neurons for related models has been studied in~\cite{praMcKeanVlasovLimitInteracting1996, delarueParticleSystemsSingular2015, fournierToyModelInteracting2015}. Following the most classical approach, we adopt the stochastic expression
\begin{equation}\label{eq:ecriture courant}
    dI(t) = bN(t)\,dt+\sqrt{2}\,dB_t.
\end{equation}

Taking the standard deviation $\sigma_c=\sqrt{2}$ constant is itself an approximation, as already
questioned in \cite{caceresAnalysisNonlinearNoisy2010}; it should in principle depend on $N(t)$. We follow most of the literature on the integrate-and-fire model in
adopting this simplification, which is more tractable.

%The sign of the mean connectivity $b$ dictates the nature of the system: It is globally excitatory when $b>0$, and globally inhibitory when $b<0$. The function $N(t)$ denotes \textbf{the firing rate of the network} also denoted \textbf{the activity function }.

Although it relies on simplifying assumptions, this expression has a major advantage: the activity function $N(t)$ is directly related to the law of the process $(V_t)_t$, being its outgoing probability flux at $+\infty$.

Combining \eqref{eq:EDS pot} and \eqref{eq:ecriture courant}, we obtain
\begin{equation}\label{eq:EDS pot2}
    \left\{\begin{aligned}
        & dV_t = \big(h(V_t)+bN(t)\big)\,dt+\sqrt{2}\,dB_t,       \\
        & \text{when} \ V_{t^-} = +\infty, \quad V_{t^+} = 0,\\
        & N(t) = \text{Flux}_{+\infty}(V_t).
      \end{aligned}\right.
\end{equation}
Setting $u(t) = \mathrm{Law}(V_t)$, we recover Eq.~\eqref{eq:FPSL Non-lin}.

To the best of our knowledge, with \cite{perthameFokkerPlanckEquationSuperlinear2026}, this is the first study of the Fokker--Planck equation associated
with the SLIF model from the point of view of partial differential equations. Apart from the
derivation of the model
itself~\cite{bretteAdaptiveExponentialIntegrateandfire2005, brunelFiringRateNoisy2003}, only the
quadratic case has been analyzed mathematically, and under the additional assumption that the
reset potential is $V_R=-\infty$ rather than $V_R=0$: imposing a Lorentzian profile on the
solutions reduces the dynamics to a coupled system of two scalar ordinary differential
equations~\cite{montbrioMacroscopicDescriptionNetworks2015}, so that the full structure of the
solutions is left aside.

\paragraph{Organization of the paper.} 
In a first section, we present our main results. Existence and uniqueness are stated in Section~\ref{sec:existenceetunicite}. Section~\ref{sec:estimates} is devoted to various estimates including large time behavior. Some technical aspects are left to appendix~\ref{appendix}.

%-------------------------------------------
\section{Notations, assumptions and main results}
\label{sec:amr}

Because the activity $N(t)$ is both defined as the flux escaping at $+\infty$ and used as the reset amplitude at $0$, Eq.~\eqref{eq:FPSL Non-lin} is conservative. It is also positivity preserving,  for all $t\geq 0$
\begin{equation}\label{eq:conservation de la masse}
 \frac{d}{dt}\int_\R u(t)\,dx=0 \qquad \text{and} \qquad u^0\geq 0 \; \Rightarrow \; u(t)\geq 0.
\end{equation}
Hence, when $u^0$ is a probability density, then, so is $u(t)$ for any $t\geq 0$. 
Apart from Section~\ref{sec:existenceetunicite}, where existence and uniqueness are established
for a general datum $u^0$ in a suitable function space, without any condition on $\int u^0$, we
shall always assume that $u$ is a probability density.

 The assumptions below are those introduced in~\cite{perthameFokkerPlanckEquationSuperlinear2026} and are satisfied by the classical superlinear integrate-and-fire drifts considered in the literature, including drifts with quadratic or exponential asymptotic growth. In particular, they ensure that the deterministic dynamics driven by $h$ can reach $+\infty$ in finite time. We assume that, for some $x_1\geq 0$,
\begin{equation}\label{eq:behavior h +infty}
    \begin{cases}
        \forall x \geq x_1, \quad h(x)>0,\quad h'(x)>0,
        \\[4pt]
        \displaystyle \int_{x_1}^{\infty}\frac{dx}{h(x)}<\infty,
        \\[8pt]
        h(x)\displaystyle\sup_{y\geq x}\frac{h'(y)^2}{h(y)^4}\in L^1(x_1,\infty).
    \end{cases}
\end{equation}
At negative potentials, we assume that $h$ has an asymptotically linear behavior,
\begin{equation}\label{eq:behavior h -infty}
    |h(x)+x|
    \underset{x\rightarrow-\infty}{=}O(1) \quad \text{and} \quad 
    \lim_{x\rightarrow-\infty}|h'(x)+1|
    =0.
\end{equation}

This hypothesis is less restrictive than the one used in \cite{perthameFokkerPlanckEquationSuperlinear2026}. It allows the class of functions $h$ to be invariant under any translation, which is convenient for establishing the existence of stationary states.

We also require the following additional growth assumption on $h$ and its derivative, 
\begin{equation}\label{eq:borne h'/h}
    \sup_{x\geq0}\frac{h'(x)}{1+|h(x)|}<\infty.
\end{equation}
Some results also use the Poincaré inequality and require that 
\begin{equation}\label{eq:croissance h Poincaré}
    \sup_{x\geq x_1}h(x)\int_x^\infty\frac{dy}{h(y)}<\infty.
\end{equation}
Assumptions \eqref{eq:behavior h +infty}, \eqref{eq:behavior h -infty} and
\eqref{eq:borne h'/h} are in force throughout the rest of the paper. Assumption
\eqref{eq:croissance h Poincaré}, which further restricts the asymptotic profile of $h$, is
recalled explicitly in every statement where it is used.

The properties of $h$ are mostly used through $u_\infty$, the stationary state of the linear problem obtained by setting $b=0$ in \eqref{eq:FPSL Non-lin}, which is a
probability density and whose properties are established
in~\cite{perthameFokkerPlanckEquationSuperlinear2026}. We only recall here the most important
one, namely its asymptotic behavior at $\pm\infty$:

\begin{equation}
    \exists \lambda_1,\lambda_2>0,  \quad \fa x \leq 0,\:\: \lambda_1e^{-x^2}\leq u_\infty(x)\leq \lambda_2e^{-\frac{x^2}{4}}\quad \text{and}\quad \exists N_\infty>0,\quad u_\infty(x)\underset{x\rightarrow\infty}{\sim}\frac{N_\infty}{h(x)}.
\end{equation}
Furthermore, it follows from \cite{perthameFokkerPlanckEquationSuperlinear2026} that when $h$ satisfies \eqref{eq:borne h'/h}, then $u_\infty >0$ satisfies
\begin{equation}\label{eq: borne u'/u}
\exists \kappa_\infty > 0, \ \fa x\in \R, \quad \big(-u_\infty'(x)\big)\leq \kappa_\infty u_\infty(x).
\end{equation}
When $h$ satisfies \eqref{eq:croissance h Poincaré}, then $u_\infty$ is a weight on $\R$ for which the following Poincaré inequality holds:
\begin{equation}\label{eq:Poincaré uinfty}
    \exists\, \alpha>0,\quad \forall\, f\in \mathcal{C}^\infty_c(\mathbb{R}),\quad
    \int_\mathbb{R} \left(f(x)-\int_\mathbb{R} f(y)\,u_\infty(y)\,dy\right)^{\!2}
    u_\infty(x)\,dx
    \leq \alpha\int_\mathbb{R} \bigl(f'(x)\bigr)^2 u_\infty(x)\,dx.
\end{equation}

When $h$ satisfies both \eqref{eq:borne h'/h} and \eqref{eq:croissance h Poincaré}, then there exist $x_2,M_1,M_2,\alpha_1,\alpha_2>0$ such that
\begin{equation}\label{eq:growth adapted h}
    \fa x \geq x_2,\quad M_1 e^{\alpha_1 x}\leq h(x) \leq M_2 e^{\alpha_2 x}.
\end{equation}

 To state our main theorems, let us first introduce the following weighted Hilbert spaces, adapted to the entropy estimates used in the proof,
\[
L^2_{\mathrm w}(\mathbb{R})=\left\{f:\mathbb{R}\rightarrow\mathbb{R}\: \Big/\:
\int_{\mathbb{R}}f(x)^2u_\infty^{-1}(x)\,dx<\infty\right\},
\]
and

\[
H^1_{\mathrm w}(\mathbb{R})=\left\{f\in L^2_{\mathrm w}(\mathbb{R})\: \Big/\:
\int_{\mathbb{R}}\Big(\partial_x\frac{f}{u_\infty}\Big)^2u_\infty\,dx<\infty\right\}.
\]

We endow $L^2_{\mathrm w}(\mathbb{R})$ with the relative entropy norm 
\begin{equation} \label{L2wNorm}
\|f\|^2_{L^2_{\mathrm w}}=\int_{\mathbb{R}}f(x)^2u_\infty^{-1}(x)\,dx .
\end{equation}

Our first main result establishes that Eq.~\eqref{eq:FPSL Non-lin} is globally well posed in
$L^2_{\mathrm w}(\R)$. 
\begin{theorem}\label{th:existence and uniqueness FPSLNL}
Let $b\in\R$ and let $u^0\in L^2_{\mathrm w}(\R)$ be a probability density. Set $w(t):=\frac{u(t)}{u_\infty}$. Then \eqref{eq:FPSL Non-lin} admits a unique global solution
\[
  (u,N)\in\Bigl(L^\infty_{\mathrm{loc}}\bigl(\R^+;L^2_{\mathrm w}(\R)\bigr)
  \cap L^2_{\mathrm{loc}}\bigl(\R^+;H^1_{\mathrm w}(\R)\bigr)\Bigr)
  \times L^2_{\mathrm{loc}}\bigl(\R^+\bigr),
\]
and, for every $T>0$,
\begin{equation}\label{eq:borne L2 uniforme N}
  \sup_{t\in[0,T]}\int_\R w(t)^2u_\infty\,dx
  \;+\;\int_0^T\!\!\int_\R\big(\partial_xw(t)\big)^2u_\infty\,dx\,dt
  \;+\;\int_0^TN(t)^2\,dt\;\leq\;\mathcal{K}\big(T,|b|,\|u^0\|_{L^2_{\mathrm w}}\big),
\end{equation}
also, the usual Sobolev inequality holds, despite the weights, 
\begin{equation}\label{eq:borne L2}
\sup_{x\in \R} |u(t,x)| \in L^3_{loc} (0,\infty).
\end{equation}
\end{theorem}

The solution is thus global for every value of the connectivity $b\in\R$: neither
$\|u(t)\|_{L^2_{\mathrm w}}$ nor the activity $N(t)$ blows up in finite time.
This is the major result of this paper, and it should be contrasted with
the NNLIF model, for which global
well-posedness in the excitatory regime is known only under a smallness condition on $b$
depending on the concentration of $u^0$ near the firing threshold.

 In a second result, we establish uniform-in-time bounds on the activity function $N(t)$ in the Cesàro sense, i.e., uniform bounds on long-time averages. The exponent in the Cesàro bounds depends on the sign of the connectivity parameter $b$, as stated in the following theorem.

\begin{theorem}[Uniform bounds on $N(t)$]\label{th:cesaro bound}
Let $(u,N)$ be the solution of \eqref{eq:FPSL Non-lin} given by
Theorem~\ref{th:existence and uniqueness FPSLNL}. The activity satisfies the following Cesàro
bounds, whose form depends on the sign of the connectivity.
\begin{enumerate}
\item[\rm(i)] \emph{Inhibitory case.} For $b\le0$, there exist $C_1(u^0)>0$ and $C_2>0$, both
independent of $b$, such that
\begin{equation}\label{eq:cesaro L2}
  \forall\,T\ge0,\qquad
  \int_0^{T}N(t)^{2}\,dt\;\le\;C_1(u^0)+C_2\,T .
\end{equation}

\item[\rm(ii)] \emph{Excitatory case.} For $b>0$, there exist $C_1'(b,u^0)>0$ and $C_2'(b)>0$
such that
\begin{equation}\label{eq:cesaro LlogL}
  \forall\,T\ge0,\qquad
  \int_0^{T}\frac{N(t)}{N_\infty}\,\ln\frac{N(t)}{N_\infty}\,dt
  \;\le\;C_1'(b,u^0)+C_2'(b)\,T .
\end{equation}
\end{enumerate}
We emphasize that the constants $C_2'(b)$ and $C_2$ do not depend on $u^0$.
\end{theorem}

% \textcolor{red}{Préciser si la borne sur $N$ dépend de la donnée initiale, quitte à prendre $T$ assez grand est-ce que cela ddépend de $b$ dans le cas inhibiteur. Pas besoin de mettre  la définition avant de Cesaro, directement écrire borne $L^2$ et $f \ln f$. Attention quand $T$ est petit, on ne peut pas dire que en. moyenne $N$ est borné car on a un $(1+T)$ et non $T$.}

The first point can be compared with~\cite{carrilloQualitativePropertiesSolutions2015}, where a similar estimate was obtained for the NNLIF model.  In the excitatory case, due to the lack of a dissipative effect associated with the sign of $b$, our approach does not provide any stronger global estimate than for $N \ln N$. However, numerical simulations strongly suggest that, for any connectivity parameter $b\in\mathbb{R}$, the activity remains uniformly bounded in time, namely $N\in L^\infty(\mathbb{R}_+)$. Establishing such a result is challenging for two reasons. On the one hand, entropy methods do not provide $L^\infty$ estimates on the activity, although a comparison
argument does give the pointwise bound \eqref{eq:est L^infini N}, valid for every $b\in\R$
provided that $u^0\in\mathcal{X}$, cf.~\eqref{eq:def X}. On the other hand, the strategy developed in~\cite{carrilloQualitativePropertiesSolutions2015} for the NNLIF model, based on global supersolutions and parabolic regularization, applies only in the case $b<0$ and $0<V_R<V_F$ and does not seem to extend to the present setting.

In our current case, due to the abstract drift function, we can no longer proceed to a change of variable that turns \eqref{eq:FPSL Non-lin} into a Stefan-like problem into which parabolic regularization is possible, as it was done in \cite{carrilloClassicalSolutionsNonlinear2011}.

The existence of at least one stationary state for every value of the connectivity parameter $b$, proved in Section~\ref{sec:stationary states}, provides a first indication that the activity should remain bounded.

\noindent Our final result concerns the asymptotic convergence towards a unique stationary state in the weak connectivity regime. We denote by $u_\infty^b$ a stationary state associated to \eqref{eq:FPSL Non-lin} with connectivity $b$ that is a probability density. In the linear case ($b=0$), convergence to equilibrium was established in~\cite{perthameFokkerPlanckEquationSuperlinear2026}. This convergence persists for sufficiently small values of $|b|$.
To prove it, we follow the strategy used for the NNLIF model in  \cite{caceresAnalysisNonlinearNoisy2010, carrilloQualitativePropertiesSolutions2015}. This consists in deriving a quadratic entropy estimate, treating the nonlinear term as a perturbation, and then using a weighted Poincaré inequality to absorb this perturbation into the diffusion before applying Grönwall's lemma. In our setting, however, the weighted Poincaré inequality associated with $u_\infty^b$ is not always available. Indeed, \eqref{eq:Poincaré uinfty} holds only under the additional assumption~\eqref{eq:croissance h Poincaré}. This assumption excludes drifts with polynomial growth at $+\infty$. Consequently, our convergence result is restricted to exponential-type drifts satisfying~\eqref{eq:growth adapted h}.

To state our result we need the stationary state $(u^b_\infty,N^b_\infty )$  for a connectivity $b$. It is proved in Section~\ref{sec:stationary states}
that there is $b_- >0$ such that for $b<b_-$, \eqref{eq:FPSL Non-lin} admits a unique normalized stationary state $(u^b_\infty,N^b_\infty )$.
Then, we measure the distance to equilibrium through a weighted norm similar to  \eqref{L2wNorm}
%$\chi^2$-divergence
\begin{equation}\label{eq:chi2}
  \|f \|_{L^2_{\mathrm w,b}}^2 \;:=\;\int_\R\frac{\bigl( f(x)\bigr)^2}{u^b_\infty(x)}\,d x
  \;=\;\int_\R\biggl(\frac{f(x)}{u^b_\infty(x)}\biggr)^{\!2}u^b_\infty(x)\,d x 
\end{equation}
for the associated Hilbert space $L^2_{\mathrm w,b}(\R)$.

\begin{theorem}[Exponential relaxation in the weakly connected regime]
\label{th:cv exp}
Assume that $h$ satisfies \eqref{eq:croissance h Poincaré}. Let $u^0\in L^2_{\mathrm w}\cap L^2_{\mathrm w,b_-}(\R)$ be a probability density. Then the solution $(u,N)$ of \eqref{eq:FPSL Non-lin} given in  Theorem~\ref{th:existence and uniqueness FPSLNL} satisfies the following:
\\[5pt]
\rm(i) \emph{Inhibitory regime.} There is $b^*_i<0$, depending only on $h$ and not on the initial datum, such that for every $b\in(b^*_i,0)$ there are
$\lambda_b>0$ and $T_b(u^0)\geq0$ for which
\begin{equation}\label{eq:cv exp N inhib}
  \sup_{t\geq T_b(u^0)}e^{\lambda_b t} \|u(t)-u^b_\infty \|^2_{L^2_{\mathrm w,b}}  \le\|u^0-u^b_\infty \|^2_{L^2_{\mathrm w,b}} \qquad \text{and} \qquad  \int_0^{\infty}e^{\lambda_b t}
  \biggl(\frac{N(t)}{N^b_\infty}-1\biggr)^{\!2}\,dt<\infty.
\end{equation}
\\[5pt]
\rm(ii) \emph{Excitatory regime.} There is $b^*_e(u^0)$ such that for every $b\in[0,\min(b^+,b^*_e(u^0)))$ there exists $\lambda>0$, depending on $b$ and on $u^0$, such that for all $t\geq 0$
\begin{equation}\label{eq:cv exp u excit}
  e^{\lambda t} \|u(t)-u^b_\infty \|^2_{L^2_{\mathrm w,b}}  \le\|u^0-u^b_\infty \|^2_{L^2_{\mathrm w,b}}
\qquad \text{and} \qquad 
  \int_0^{\infty}e^{\lambda t}
  \biggl(\frac{N(t)}{N^b_\infty}-1\biggr)^{\!2}\,dt <\infty .
\end{equation}
\end{theorem}

The proofs and statements in the excitatory and inhibitory regimes are different because of the asymmetric decays at $\pm \infty$ of the weights. Especially that according to \eqref{eq: borne u'/u} only the right-tail of $|u_\infty'|$ is dominated by $u_\infty$. 
The fact that $b_i^*$ is independent of $u^0$ can be compared with the results in~\cite{carrilloQualitativeProperties2015}.

\section[Wellposedness of solutions]{Existence and uniqueness of solutions of \eqref{eq:FPSL Non-lin}}
\label{sec:existenceetunicite}

The main ingredient to prove Theorem~\ref{th:existence and uniqueness FPSLNL} is a suitable relative entropy associated with
the steady state $u_\infty$ of the linear equation corresponding to the case $b=0$.

Let $(u,N)$ satisfy, for some given input $J(t)$ and some right-hand side $R$, the equation
\begin{equation}\label{eq:generique}
	\left\{\begin{aligned}
		& \frac{\partial u}{\partial t}(t,x)
		+ \frac{\partial}{\partial x}\big((h(x)+J(t)) u(t,x)\big)
		- \frac{\partial^2 u}{\partial x^2}(t,x)
		= \delta_0(x)\,N(t)
		+R(t,x),
		\quad x \in \mathbb{R},\ t \ge 0,
		\\
		& N(t) = \lim_{x \to +\infty} h(x)u(t,x).
	\end{aligned}\right.
\end{equation}
Then an immediate application of the chain rule gives, for every smooth function $H$ and
$w=\frac{u}{u_\infty}$,
\begin{align}
\frac{d}{dt} \int_\R u_\infty H(w)\, dx =&  -\int_\R u_\infty H''(w) (\partial_x w)^2dx \notag
 -N_\infty\Big(H\big(\tfrac{N}{N_\infty}\big)-H(w(0))-H'(w(0))(w(0)-\tfrac{N}{N_\infty})\Big) \\
&-J(t)\int_\R \big(wH'(w)-H(w)\big)u_\infty'\,dx + \int_\R R(t,x) H'(w)\, dx.
\label{RE:generique}
\end{align}
The interest of this computation is that the first two terms on the right-hand side are non-positive when $H$ is convex.

The proof of Theorem~\ref{th:existence and uniqueness FPSLNL} is organized as follows. We first
consider, for a given input function $J$, the linear counterpart of \eqref{eq:FPSL Non-lin},
obtained by replacing the nonlinear drift term $bN(t)$ by the prescribed function $J(t)$. The
well-posedness of this linear problem is stated in Section~\ref{sec:linear} and proved in
Appendix~\ref{appendix}.

We then prove, in Section~\ref{sec:unicite}, the uniqueness of two global solutions with a given
initial datum $u^0\in L^2_{\mathrm w}(\R)$. Truncated quadratic entropy estimates, written in the
form of a Lipschitz inequality, produce an absorption of the error
\[
\int_0^{t_0}\big(N_1(s)-N_2(s)\big)^2\,ds
\]
on a small time interval $(0,t_0)$, which is then propagated to the whole interval of existence.

In Section~\ref{sec:existence X} we show that, for the more restrictive class of initial data
\begin{equation}\label{eq:def X}
\mathcal{X}:=\Bigl\{u^0\in L^2_{\mathrm w}(\R)\ :\ \frac{u^0}{Q}\in L^\infty(\R)\Bigr\},
\end{equation}
where $Q$ is the comparison function defined in \eqref{eq:def Q}, the Lipschitz estimate obtained
in the uniqueness proof can be improved into a uniform one. This turns the map $J\mapsto bN$ into
a genuine contraction and yields the global existence of the solution on $[0,\infty)$.

The extension of global existence to an arbitrary initial datum $u^0\in L^2_{\mathrm w}(\R)$
is obtained by a weak limit, using the usual  compactness argument of Aubin--Lions--Simon
type cf.~\cite{simonCompactSetsSpace1987}. This last step requires a uniform bound on the activities of the approximating sequence,
which is only available once the estimates of Section~\ref{sec:estimates} have been established;
it is therefore postponed to Section~\ref{sec:existence L2w suite}.

\subsection{The linear problem with a time-dependent drift $J(t)$}
\label{sec:linear}

Given an input function $J \in L^2_{\mathrm{loc}}(\R^+)$, we consider the linear equation
\begin{equation}\label{eq:FPSL with J(t)}
	\left\{\begin{aligned}
		& \frac{\partial u}{\partial t}(t,x)
		+ \frac{\partial}{\partial x}\big(\big(h(x)+J(t)\big)u(t,x)\big)
		- \frac{\partial^2 u}{\partial x^2}(t,x)
		= \delta_0(x)\,N(t),
		\quad x \in \mathbb{R},\ t \ge 0, \\
		& N(t) = \lim_{x \to +\infty} h(x)u(t,x).
	\end{aligned}\right.
\end{equation}
The weak formulation of \eqref{eq:FPSL with J(t)} is defined in the spirit of \cite{perthameFokkerPlanckEquationSuperlinear2026}, for all test functions
$\phi \in \mathcal{C}_c^\infty(\R^+\times\R)$, by the equality
\begin{align}
	\int_0^\infty \! \int_{\mathbb{R}} u(t,x)
	\Big(- \partial_t \phi(t,x)
		- \big(h(x)+& J(t)\big)\partial_x \phi(t,x)
		- \partial_x^2 \phi(t,x)\Big)\,dx\,dt \notag
\\
	&= \int_0^\infty N(t)\phi(t,0)\,dt
	+ \int_{\mathbb{R}} u^0(x)\phi(0,x)\,dx.
	\label{eq:weak formulation FPS}
\end{align}
and 
\[
\text{for a.e.}\:\:t\geq 0, \quad \int_\R u(t,x)\,dx = \int_\R u^0(x)\,dx.
\]
In order to state the estimates on solutions of \eqref{eq:FPSL with J(t)}, we introduce the
following functions:
\[
\tilde{H}(x) = \int_0^x\big(h(y)-h(0)\big)\,dy, \qquad
c_0=\Big(\int_0^\infty\exp\big(-\tilde{H}(y)\big)\,dy\Big)^{-1},
\]
\begin{equation}\label{eq:def Q}
    Q(x) =\begin{cases}
        1 &\text{for}\quad x\leq 0,
        \\[4pt]
        c_0 \exp\big(\tilde{H}(x)\big)\displaystyle\int_x^\infty\exp\big(-\tilde{H}(y)\big)\,dy
        &\text{for}\quad x>0.
    \end{cases}
\end{equation}
Before stating the well-posedness theorem for \eqref{eq:FPSL with J(t)}, we record two properties of $Q$ which we use repeatedly.

\begin{lemma}\label{lem:proprietes Q}
The function $Q$ defined by \eqref{eq:def Q} satisfies
\[
    Q(x)\underset{x\rightarrow+\infty}{\sim}\frac{c_0}{h(x)}
    \qquad\text{and}\qquad
    q_0:=\sup_{x\in\R}\Big|\frac{Q'(x)}{Q(x)}\Big|<\infty .
\]
In particular $\frac{Q^2}{u_\infty}$ is integrable in a neighborhood of $+\infty$.
\end{lemma}

\begin{theorem}[Well-posedness for the linear problem]\label{th:existence unicite avec J(t)}
Let $u^0\in L^2_{\mathrm w}(\R)$ and $J\in L^2_{\mathrm{loc}}(\R^+)$. Then
\eqref{eq:FPSL with J(t)} admits a unique weak solution
\[
  (u,N)\in L^\infty_{\mathrm{loc}}\big(\R^+;L^2_{\mathrm w}(\R)\big)\times L^2_{\mathrm{loc}}(\R^+).
\]
Furthermore, defining $w=\frac{u}{u_\infty}$ and $\mathcal{J}(T)=\int_0^{T}J(t)^2\,dt$, there is a constant
$C>0$ depending only on $h$ such that, for every $T>0$,
\begin{equation}\label{eq:a priori bound solution FPSL}
  \int_\R w(T,x)^2u_\infty(x)\,dx
  +\int_0^{T}\!\!\int_\R\big(\partial_xw(t,x)\big)^2u_\infty(x)\,dx\,dt
  \;\leq\;e^{\mathcal{J}(T)}\,\|u^0\|_{L^2_{\mathrm w}}^2,
\end{equation}
\begin{equation}\label{est: borne sur N}
  \int_0^{T}N(t)^2\,dt\;\leq\;C\big(T+e^{\mathcal{J}(T)}\big)\,\|u^0\|_{L^2_{\mathrm w}}^2,
\end{equation}
and
\begin{equation}\label{est:borne L^1 J}
    \fa \varepsilon >0,\ \exists\, \gamma_\varepsilon> 0,\quad \fa T\geq 0,\qquad
    \int_0^T|N(t)|\,dt\;\leq\;\big[ \varepsilon\int_0^T |J(t)|\,dt+\gamma_\varepsilon T+1\big] \|u^0\|_{L^1(\R)} .
\end{equation}
When, additionally, $\frac{u^0}{Q}\in L^\infty(\R)$, then the solution satisfies the $L^\infty$
bound
\begin{equation}\label{est:borne L^infini J}
    \fa t\geq0,\qquad
    \Big\|\frac{u(t)}{Q}\Big\|_{\infty}\;\leq\;\Big\|\frac{u^0}{Q}\Big\|_{\infty}\,e^{\Lambda_J(t)},
    \qquad
    \Lambda_J(t) := C\Big(t+\int_0^t|J(s)|\,ds\Big).
\end{equation}
\end{theorem}

The dependency on the initial data in the $L^1$ norm only is fundamental for the nonlinear problem. 
\begin{remark}[The nonlinearity $J=bN$]\label{rem:bornes nonlin}
Both estimates of Theorem~\ref{th:existence unicite avec J(t)} can be used for solutions of the nonlinear problem \eqref{eq:FPSL Non-lin}, that is with $J=bN$, when $u^0$ is a probability density. Choosing
$\varepsilon=\frac{1}{2|b|}$ in \eqref{est:borne L^1 J} gives
\begin{equation}\label{eq:est L^1 N}
    \fa T\geq0,\qquad \int_0^T N(t)\,dt\;\leq\;2\big(\gamma_bT+1\big),
    \qquad \gamma_b:=\gamma_{\frac{1}{2|b|}} ,
\end{equation}
and \eqref{est:borne L^infini J} then yields the pointwise bound
\begin{equation}\label{eq:est L^infini N}
    \fa t\geq0,\qquad N(t)\;\leq\;C^0_b\,e^{a_bt},
    \qquad C^0_b:=c_0\Big\|\frac{u^0}{Q}\Big\|_\infty e^{2C|b|},
    \qquad a_b:=C\big(1+2|b|\gamma_b\big) ,
\end{equation}
for every $b\in\R$. Notice that \eqref{eq:est L^1 N} holds for
any probability density $u^0\in L^2_{\mathrm w}(\R)$, whereas
\eqref{eq:est L^infini N} requires in addition $u^0\in\mathcal{X}$.
\end{remark}

We only prove here the two comparison estimates \eqref{est:borne L^infini J} and
\eqref{est:borne L^1 J}. The remaining bounds, together with the construction of the solution, are established in
Appendix~\ref{appendix} using a piecewise constant approximation of $J(t)$.

\begin{proof}[Proof of \eqref{est:borne L^infini J} and \eqref{est:borne L^1 J}]
\emph{The $L^\infty$ bound.} Set $A:=\|h'\|_{\infty,(-\infty,0)}$, and recall
$q_0=\|Q'/Q\|_\infty$ from Lemma~\ref{lem:proprietes Q}. Let
\[
  \lambda(t):=\max\Big(A,\ \big|h(0)+J(t)\big|\,q_0\Big),
  \qquad \Lambda(t):=\int_0^t\lambda(\tau)\,d\tau,
  \qquad P(t,x):=e^{\Lambda(t)}Q(x).
\]
Notice that $\lambda(t)\leq A+|h(0)|q_0+q_0|J(t)|$, so that
$\Lambda(t)\leq C\big(t+\int_0^t|J(s)|ds\big)= \Lambda_J(t)$ with $C:=\max\big(A+|h(0)|q_0,\,q_0\big)$.
\\

We claim that $P$ is a supersolution. For $x<0$ we have $Q\equiv1$, hence
\[
  \partial_t P+\partial_x\big((h+J(t))P\big)-\partial_x^2P
  =e^{\Lambda(t)}\big(\lambda(t)+h'(x)\big)\;\geq\;e^{\Lambda(t)}\big(\lambda(t)-A\big)\;\geq\;0 .
\]
For $x>0$, the definition \eqref{eq:def Q} gives $\big(h(x)-h(0)\big)Q(x)-Q'(x)=c_0$, so that
$(hQ)'-Q''=h(0)Q'$ and therefore
\[
  \partial_t P+\partial_x\big((h+J(t))P\big)-\partial_x^2P
  =e^{\Lambda(t)}Q\Big(\lambda(t)+\big(h(0)+J(t)\big)\frac{Q'}{Q}\Big)
  \geq e^{\Lambda(t)}Q\Big(\lambda(t)-\big|h(0)+J(t)\big|q_0\Big)\geq0 .
\]
Moreover $P(t,\cdot)$ is continuous, with $\partial_xP(t,0^-)=0$ and
$\partial_xP(t,0^+)=-c_0e^{\Lambda(t)}$, while $\lim_{x\to\infty}h(x)P(t,x)=c_0e^{\Lambda(t)}$;
hence $\lim_{x\to\infty}h(x)P(t,x)=\partial_xP(t,0^-)-\partial_xP(t,0^+)$ and $P$ satisfies
\[
  \partial_t P+\partial_x\big((h+J(t))P\big)-\partial_x^2P\;\geq\;N_P(t)\,\delta_0(x),
  \qquad N_P(t)=\lim_{x\to\infty}h(x)P(t,x).
\]
Defining $L_Jf := -\partial_x\big((h+J(t))f\big)+f''+\big(\lim_{x\to\infty}h(x)f(x)\big)\delta_0$,
the functions $u$ and $P$ therefore satisfy 
\[
\partial_tu-L_Ju=0 \qquad \text{and} \qquad  \partial_tP-L_JP=g\geq 0.
\]
With $C_0:=\big\|\tfrac{u^0}{Q}\big\|_\infty$, it holds $u^0\leq C_0 P(0,\cdot)$.
From this it follows that $u \leq C_0 P$ by the maximum principle and \eqref{est:borne L^infini J} follows. 

For the sake of completeness, we prove the maximum principle in our context. This method is inspired by the theory of entropy for parabolic equations that can be found in \cite{michelGeneralRelativeEntropy2005} and \cite{perthameTransportEquationsBiology2007}. Choose the entropy 
$H(x):=(x-C_0)_+^2$, set $\omega:=\tfrac{u}{P}$, and use the relative entropy of $u$ with respect to $P$:
\[
\partial_t\big(PH(\omega)\big)=L_JP\big(H(\omega)-\omega H'(\omega)\big)+L_Ju\, H'(\omega)
+g\big(H(\omega)-\omega H'(\omega)\big).
\]
With this choice of $H$ one has $H(\omega)-\omega H'(\omega)=-(\omega+C_0)(\omega-C_0)_+\leq0$,
so that the last term is non-positive since $g\geq0$. Moreover, writing
$\omega_\infty=\lim_{x\to\infty}\omega$, a direct computation gives
\begin{align*}
        \int_\R \Big(L_JP\big(H(\omega)- & \omega H'(\omega)\big)+L_Ju\, H'(\omega) \Big)dx\\
    &= - \int_\R\big(\partial_x\omega\big)^2H''(\omega)P\,dx
    -N_P(t)\Big(H(\omega_\infty)-H(\omega(0))-H'(\omega(0))\big(\omega_\infty-\omega(0)\big)\Big)
    \leq 0 .
\end{align*}
Hence $t\mapsto\int_\R H(\omega)P\,dx$ is non-increasing, and since it vanishes at $t=0$ by
definition of $C_0$, it vanishes for every $t\geq0$. Therefore $u(t)\leq C_0P(t)$, that is the maximum principle.

\medskip
\emph{The $L^1$ bound.} Let  $L\geq 2$. Let
$\chi\in\Cc^\infty(\R)$ be a positive increasing smooth function that satisfies
\[
0 \leq \chi'\leq \frac{2}{L}, \quad \chi = 0 \text{ on } (-\infty,0], \quad
\chi(x) = \frac{x}{L} \text{ on } [1,L-1], \quad  \chi = 1 \text{ on } [L,\infty).
\]
The pair $(|u|,|N|)$ satisfies, in the distributional sense,
\[
\partial_t |u| + \partial_x \Big( \big( J(t) + h(x) \big) |u| \Big) - \partial_x^2 |u|
\leq |N(t)|\, \delta_0 .
\]
Multiplying by $\chi$, and using $\chi(0)=0$ together with
$\lim_{x\to\infty}h(x)|u(t,x)|=|N(t)|$, integration by parts  over $\R$ gives
\[
    \frac{d}{dt}\int_\R|u(t,x)|\chi(x)\,dx+|N(t)|
    \leq \int_\R\big((J(t)+h(x))\chi'(x)+\chi''(x)\big)|u(t,x)|\,dx
    \leq \big[ \|h\chi'+\chi''\|_\infty+\frac{|J(t)|}{L}\big] \|u^0\|_{L^1(\R)},
\]
where we used $\int_\R|u(t,x)|dx\leq \|u^0\|_{L^1(\R)}$. We now integrate
over $(0,T)$, use $0\leq\chi\leq1$ together with $|u(T)|\chi\geq0$ and for $\eps >0$ we set $L=\max(2,\frac 1\eps)$. We obtain
\[
    \int_0^T|N(t)|\,dt\;\leq\;\int_\R|u^0|\chi\,dx+\big[\gamma_\varepsilon T+\varepsilon\int_0^T|J(t)|\,dt\big] \|u^0\|_{L^1(\R)}
    \;\leq\;\big[ 1+\gamma_\varepsilon T+ \varepsilon\int_0^T|J(t)|\,dt \big] \|u^0\|_{L^1(\R)},
\]
with $\gamma_\varepsilon:=\|h\chi'+\chi''\|_\infty$, which is \eqref{est:borne L^1 J}. 

\end{proof}

\subsection{Uniqueness of solutions of \eqref{eq:FPSL Non-lin}}
\label{sec:unicite}

In order to prove uniqueness, we first establish a continuity estimate between two solutions of
the linear problem associated with two different drifts.

\begin{prop}\label{prop:controle I}
Let $u^0\in L^2_{\mathrm w}(\R)$ and $J_1,J_2\in L^2_{\mathrm{loc}}(\R^+)$, and let $(u_1,N_1)$
and $(u_2,N_2)$ be the solutions of \eqref{eq:FPSL with J(t)} associated respectively with the
drift functions $J_1$ and $J_2$, with $u_i(0) = u^0$ for $i=1,2$. Then, 
$w_i(t)=\frac{u_i(t)}{u_\infty}$ satisfy for any $t\geq 0$,
\begin{equation}\label{eq: controle I}
    \int_{\R}\big(w_{1}(t)-w_{2}(t)\big)^{2} u_{\infty}\,dx
    \;\leq\; \|u^0\|_{L^2_{\mathrm w}}^{2}\;
    e^{\mathcal{J}_1(t)+\mathcal{J}_2(t)}
    \int_{0}^{t}\big(J_{1}(\tau)-J_{2}(\tau)\big)^{2}\,d\tau .
\end{equation}
\end{prop}

\begin{proof}
The difference $(u,N)=(u_1-u_2, N_1-N_2)$ satisfies \eqref{eq:generique} with
\[
R(t,x) = J_2(t) \partial_x u_2- J_1(t) \partial_x u_1
=-\big(J_1(t)-J_2(t)\big) \partial_x u_2 -J_1(t) \partial_x (w u_\infty).
\]
We set $w(t)=w_1-w_2$ and $I(t) = \int_{\R}w(t)^2 u_\infty \, dx$, and we use the relative
entropy identity \eqref{RE:generique} with $H(w)=w^2$. After integration by parts in the last
two terms, we find
\begin{align*}
    \frac{d}{dt}I(t) &=  -2 \int_{\R} \big(\partial_x w\big)^2 u_\infty \, dx
    - N_\infty \Big( \tfrac{N}{N_\infty}- w(t,0)\Big)^2
     +2\big(J_1(t)-J_2(t)\big) \int_{\R} \partial_x w \;  w_2\, u_\infty\,dx
     + 2 J_1(t) \int_{\R} \partial_x w  \; w\, u_\infty\,dx .
\end{align*}
It remains to estimate the last two terms. By the Cauchy--Schwarz inequality and
$ab \leq \tfrac{a^2}{2}  + \tfrac{b^2}{2}$,
\[
    2\big|J_1 - J_2\big|\, \Big| \int_{\R} \partial_x w \; w_2 \,u_\infty \, dx \Big|
    \leq \int_{\R} \big(\partial_x w\big)^2 u_\infty \, dx
    + \big(J_1 - J_2\big)^2 \int_{\R} u_\infty w_2^2 \, dx ,
\]
and, in the same way,
\[
    2\big|J_1\big| \,\Big|\int_{\R} \partial_x w \; w \,u_\infty \, dx\Big|
    \leq  \int_{\R} \big(\partial_x w\big)^2 u_\infty \, dx
    + J_1^2\, I(t) .
\]
Both diffusion contributions are absorbed by $-2\int_\R(\partial_xw)^2u_\infty\,dx$, and
discarding the remaining non-positive terms we obtain
\[
I'(t) \leq \big(J_1(t) - J_2(t)\big)^2 \|u_2(t)\|_{L^2_{\mathrm w}}^{2}
+ J_1(t)^2\, I(t) .
\]
By \eqref{eq:a priori bound solution FPSL},  we know that $\|u_2(t)\|_{L^2_{\mathrm w}}^{2}
\leq e^{\mathcal{J}_2(t)}\|u^0\|_{L^2_{\mathrm w}}^{2}$, so that
\[
I'(t) - J_1(t)^2 I(t)
\leq \big(J_1(t) - J_2(t)\big)^2 e^{\mathcal{J}_2(t)}
\|u^0\|_{L^2_{\mathrm w}}^{2}.
\]
Multiplying by the integrating factor $e^{-\mathcal{J}_1(t)}$, integrating and using
$I(0)=0$, we get
\[
I(t) \leq \|u^0\|_{L^2_{\mathrm w}}^{2}\,
e^{\mathcal{J}_1(t)+\mathcal{J}_2(t)}
\int_0^t \big(J_1(\tau) - J_2(\tau)\big)^2 \, d\tau ,
\]
which implies \eqref{eq: controle I}.
\end{proof}

Uniqueness is now obtained by showing that, on a small time interval, the quantity
$\int_0^{t_0}|N_1-N_2|^2$ is bounded by a fraction of itself, hence vanishes. The key point is
that the two solutions are \emph{given}: no uniformity with respect to a class of data is
required, and the smallness comes from the absolute continuity of the Lebesgue integral.

\begin{prop}[Uniqueness for \eqref{eq:FPSL Non-lin}]\label{prop:unicite L2w}
Let $b\in\R$. Two solutions of \eqref{eq:FPSL Non-lin} on a
common interval, with the same initial datum $u^0\in L^2_{\mathrm w}(\R)$, both belonging to the class of
Theorem~\ref{th:existence and uniqueness FPSLNL}, coincide.
\end{prop}

\begin{proof}
For $(u_1,N_1)$ and $(u_2,N_2)$ two such solutions, we first show that $(u_1,N_1)=(u_2,N_2)$ on a small interval $[0,t_0]$.

Let $\varphi$ satisfy
\begin{equation}\label{def:varphi}
\varphi\in\Cc^\infty(\R, \R^+), \qquad \varphi' > 0 \text{ in } (x_1,x_1+1), \qquad
\varphi|_{(-\infty,x_1)}=0 \quad \text{and} \quad \varphi|_{(x_1+1,\infty)}=1 .
\end{equation}
Taking supp$(\varphi)\subset [x_1,\infty)$ ensures that $\varphi u_\infty'\leq 0$. 
We then set
\[
I_\varphi(t) = \int_{\R}\big(w_{1}(t)-w_{2}(t)\big)^{2} u_{\infty}\varphi\,dx,
\qquad
I(t) = \int_{\R}\big(w_{1}(t)-w_{2}(t)\big)^{2} u_{\infty}\,dx .
\]
Applying \eqref{RE:generique} with $H(x)=x^2$, the truncation $\varphi$ and
$R=-b\big(N_1\partial_xu_1-N_2\partial_xu_2\big)$ gives the {\em truncated quadratic entropy}
\begin{equation}\label{eq:derivation Iphi}
\begin{aligned}
    I_\varphi'(t) ={} &-2\int_\R \big(\partial_x(w_1-w_2)\big)^2u_\infty\varphi \,dx
    +N_\infty\bigg(-\Big(\frac{N_1-N_2}{N_\infty}\Big)^2+\int_0^\infty (w_1-w_2)^2\,\varphi'\,dx\bigg)\\
    &+\int_\R (w_1-w_2)^2\big(u_\infty \varphi '\big)'\,dx
    -2b(N_1-N_2)\int_\R (w_1-w_2)\,\partial_x(u_\infty w_1)\,\varphi  \,dx\\
    &-2bN_2\int_\R (w_1-w_2)\,\partial_x\big(u_\infty (w_1-w_2)\big)\,\varphi  \,dx.
\end{aligned}
\end{equation}
Setting $\theta_1:=\big\|\tfrac{\varphi'}{u_\infty}\big\|_\infty$ and
$\theta_2:=\big\|\tfrac{u'_\infty \varphi'}{u_\infty}+\varphi''\big\|_\infty$, we have
$\varphi'\leq\theta_1u_\infty$ and $\big(u_\infty\varphi'\big)'\leq \theta_2 u_\infty$, hence
\[
  N_\infty\int_0^\infty (w_1-w_2)^2\varphi'\,dx\leq N_\infty\theta_1\,I(t),
  \qquad
  \int_\R (w_1-w_2)^2\big(u_\infty\varphi'\big)'\,dx\leq \theta_2\,I(t).
\]
We next treat the fourth term of \eqref{eq:derivation Iphi} by the Cauchy--Schwarz inequality
and $(a+b)^2\leq 2(a^2+b^2)$, using $|u_\infty'|\leq\kappa_\infty u_\infty$ on
$\operatorname{supp}\varphi\subset[x_1,\infty)$,
\begin{align*}
       -2b(N_1-N_2)\int_\R (w_1-w_2)\,\partial_x(u_\infty w_1)\,\varphi\,dx
    &\leq 2\big|b(N_1-N_2)\big|\int_\R |w_1-w_2|\big(\kappa_\infty w_1 +|\partial_xw_1|\big)u_\infty\,dx\\
    &\leq 2\big|b(N_1-N_2)\big|\,I(t)^{1/2}\bigg(\int_\R\big(\kappa_\infty w_1+|\partial_x w_1|\big)^2u_\infty\,dx\bigg)^{1/2}\\
    &\leq 2\sqrt2\,\big|b(N_1-N_2)\big|\,I(t)^{1/2}
    \bigg(\int_\R\big((\kappa_\infty w_1)^2+(\partial_x w_1)^2\big)u_\infty\,dx\bigg)^{1/2}.
\end{align*}
Setting
\[
F(w_1(t)) := 8b^2\int_\R\big((\kappa_\infty w_1(t))^2+(\partial_x w_1(t))^2\big)u_\infty\,dx,
\]
the right-hand side above is exactly $\big|N_1(t)-N_2(t)\big|\,I(t)^{1/2}F(w_1(t))^{1/2}$. Since
$u_1$ belongs to $L^\infty_{\mathrm{loc}}(L^2_{\mathrm w})\cap L^2_{\mathrm{loc}}(H^1_{\mathrm w})$,
the function $t\mapsto F(w_1(t))$ belongs to $L^1_{\mathrm{loc}}(\R^+)$, and therefore, by
absolute continuity of the Lebesgue integral,
\begin{equation}\label{eq: F dans L1loc}
    \lim_{T\rightarrow0}\int_0^T F(w_1(t))\,dt =0.
\end{equation}
Finally, integrating by parts in the last term of \eqref{eq:derivation Iphi},
\begin{align*}
    -2bN_2\int_\R (w_1-w_2)\,\partial_x\big(u_\infty (w_1-w_2)\big)\varphi \,dx
    &=-bN_2\int_\R (w_1-w_2)^2\big(u_\infty'\varphi-u_\infty \varphi'\big)\,dx\\
    &\leq |N_2(t)|\,\underbrace{|b|\big(\kappa_\infty+\|\varphi\|_\infty\big)}_{:=\theta_3}
    \int_\R (w_1-w_2)^2u_\infty\,dx
    = \theta_3\,|N_2(t)|\,I(t).
\end{align*}
Gathering all the previous estimates and discarding the non-positive diffusion term, we obtain
\[
    I'_\varphi(t)+\frac{\big(N_1(t)-N_2(t)\big)^2}{N_\infty}
    \leq \big(N_\infty\theta_1+\theta_2+\theta_3|N_2(t)|\big)I(t)
    +\big|N_1(t)-N_2(t)\big|\,I(t)^{1/2}F(w_1(t))^{1/2}.
\]
Consider now a time $T>0$, integrate over $(0,T)$ and use $I_\varphi(0)=0$ and $I_\varphi(T)\geq 0$,
\[
    \int_0^T\frac{\big(N_1(t)-N_2(t)\big)^2}{N_\infty}\,dt
    \leq \int_0^T\big(N_\infty\theta_1+\theta_2+\theta_3|N_2(t)|\big)I(t)\,dt
    +\int_0^T\big|N_1(t)-N_2(t)\big|\,I(t)^{1/2}F(w_1(t))^{1/2}\,dt .
\]
By Proposition~\ref{prop:controle I} applied with $J_i=bN_i$, we have, for every $t\in(0,T)$,
\[
I(t)\leq M(T)\int_0^T\big(N_1(s)-N_2(s)\big)^2\,ds,
\qquad
M(T):= b^2\|u^0\|_{L^2_{\mathrm w}}^2\exp\bigg(b^2\int_0^T\big(N_1^2+N_2^2\big)\,ds\bigg).
\]
Two applications of the Cauchy--Schwarz inequality then give
\begin{align*}
    \int_0^T\frac{\big(N_1(t)-N_2(t)\big)^2}{N_\infty}\,dt
    &\leq M(T)\int_0^T\big(N_\infty\theta_1+\theta_2+\theta_3|N_2(t)|\big)\,dt
    \int_0^T\big(N_1-N_2\big)^2\,ds\\
    &\quad+M(T)^{1/2}\bigg(\int_0^TF(w_1(t))\,dt\bigg)^{1/2}\int_0^T\big(N_1-N_2\big)^2\,ds ,
\end{align*}
where we used $\int_0^T|N_1-N_2|F(w_1)^{1/2}dt\leq\big(\int_0^T(N_1-N_2)^2\big)^{1/2}
\big(\int_0^TF(w_1)\big)^{1/2}$. Setting
\[
\varepsilon(T):=M(T)\Big((N_\infty\theta_1+\theta_2)T
+\theta_3\sqrt{T}\Big(\int_0^TN_2^2\,dt\Big)^{1/2}\Big)
+M(T)^{1/2}\Big(\int_0^TF(w_1(t))\,dt\Big)^{1/2},
\]
the continuity of $T\mapsto\int_0^TN_i(t)^2\,dt$ for $i=1,2$ and \eqref{eq: F dans L1loc} ensure
that $\lim_{T\rightarrow 0}\varepsilon(T)=0$. Choosing $t_0>0$ such that
$N_\infty\varepsilon(t_0)<1$, the inequality above reads
$X\leq N_\infty\varepsilon(t_0)X$ with $X=\int_0^{t_0}(N_1-N_2)^2$, whence
\[
\int_0^{t_0}\big(N_1(s)-N_2(s)\big)^2\,ds=0 .
\]
Therefore $N_1=N_2$ almost everywhere on $(0,t_0)$ and, by the uniqueness part of Theorem~\ref{th:existence unicite avec J(t)}, $u_1=u_2$ on $(0,t_0)$.

To conclude, let  
$\tau=\sup \bigl\{t\geq 0\ :\ u_1|_{(0,t)}=u_2|_{(0,t)}\bigr\}\geq t_0>0$ and assume $\tau<\infty$. Since $N_1=N_2$ almost everywhere on
$[0,\tau)$ and $\mathcal{J}_i(\tau)<\infty$, the right-hand side of \eqref{eq: controle I}
vanishes at $t=\tau$, so that $u_1(\tau)=u_2(\tau)$ in $L^2_{\mathrm w}(\R)$. We may repeat the
above argument with this common value as initial datum and propagates the equality strictly beyond
$\tau$, a contradiction. Hence $\tau=\infty$, which concludes the proof.
\end{proof}

\subsection{Global existence for an initial datum in $\mathcal{X}$} 
\label{sec:existence X}

We have gathered the material to prove Theorem~\ref{th:existence and uniqueness FPSLNL} for an initial datum
$u^0\in\mathcal{X}$. Taking the datum in this class allows us to use the $L^\infty$ estimate
\eqref{est:borne L^infini J}, which, in contrast with the estimate of the previous proof,
produces a Lipschitz constant that is \emph{uniform} on balls of $L^2(0,T)$. This is the content
of the following proposition, whose proof rests on a cutoff function pushed towards $+\infty$.

\begin{prop}\label{prop:linear_locLipschitz}
Let $R>0$, $T\in (0,1)$ and let $u^0\in \mathcal{X}$ such that $u^0\geq 0$. There is a constant $C(R,T,u^0)>0$ such that, for any $J_1,J_2 \in L^2(0,T)$ with $\|J_1\|_{2}\leq R$, $\|J_2\|_{2}\leq R$, we have
\begin{equation}\label{eq:N-J lipschitz}
    N_\infty\int_0^T\Big(\frac{N_1-N_2}{N_\infty}\Big)^2\,dt\leq C(R,T,u^0)\,\|J_1-J_2\|_2^2, \qquad \quad \lim_{T\rightarrow 0}C(R,T,u^0)=0.
\end{equation}
The constant $C(R,T,u^0)$ depends on $u^0$ only through $\big \|\tfrac{u^0}{Q}\big\|_\infty$ and $\|u^0\|_{L^2_{\mathrm w}}$.
\end{prop}

The major difficulty compared to the proof of Proposition~\ref{prop:unicite L2w} is that we cannot rely on time compactness. The control on the constant $C(R,T,u^0)$ is more demanding, which requires a more involved argument.

\begin{proof}
Let $J_1,J_2 \in L^2(0,T)$ with $\max_{i=1,2}\|J_i\|_{2} \leq R$, and let $(u_i,N_i)$ be the
solution of \eqref{eq:FPSL with J(t)} associated with the drift $J_i$. As in the proof of
Proposition~\ref{prop:unicite L2w} we use a truncated quadratic entropy, but in order to exploit
the condition $u^0\in\mathcal{X}$ we introduce a parameter $X>0$ and the shifted cutoff
$\varphi_X(x) := \varphi(x-X)$, where $\varphi$ is given by \eqref{def:varphi}. Since
$\varphi\geq0$ and $\varphi''$ is bounded, one has $(\varphi')^2\leq2\|\varphi''\|_\infty\varphi$,
so that
\[
\Phi_0:=\Big\|\varphi+\tfrac{(\varphi')^2}{\varphi}\Big\|_\infty<\infty .
\]
We write $I_X(t)=\int_\R\big(w_1(t)-w_2(t)\big)^2u_\infty\varphi_X\,dx$, and $I(t)$, $w_i$ are as
before. Replacing $bN_i$ by $J_i$, the identity \eqref{eq:derivation Iphi} becomes
\begin{equation}\label{eq:derivation IphiX}
\begin{aligned}
    I_X'(t) ={} &-2\int_\R \big(\partial_x(w_1-w_2)\big)^2u_\infty\varphi_X \,dx
    +N_\infty\bigg(-\Big(\frac{N_1-N_2}{N_\infty}\Big)^2+\int_0^\infty (w_1-w_2)^2\,\varphi_X'\,dx\bigg)\\
    &+\int_\R (w_1-w_2)^2\big(u_\infty\varphi_X '\big)'\,dx
    -2(J_1-J_2)\int_\R (w_1-w_2)\,\partial_x(u_\infty w_1)\,\varphi_X  \,dx\\
    &-2J_2\int_\R (w_1-w_2)\,\partial_x\big(u_\infty (w_1-w_2)\big)\,\varphi_X  \,dx.
\end{aligned}
\end{equation}
Setting $\theta_1(X):=\big\|\tfrac{\varphi'_X}{u_\infty}\big\|_\infty$ and
$\theta_2(X):=\big\|\tfrac{u'_\infty \varphi'_X}{u_\infty}+\varphi''_X\big\|_\infty$, we get
$\varphi_X'\leq\theta_1(X)u_\infty$ and $\big(u_\infty\varphi'_X\big)'\leq \theta_2(X) u_\infty$,
hence
\[
  N_\infty\int_0^\infty (w_1-w_2)^2\varphi_X'\,dx\leq N_\infty\theta_1(X)\,I(t),
  \qquad
  \int_\R (w_1-w_2)^2\big(u_\infty\varphi_X'\big)'\,dx\leq \theta_2(X)\,I(t).
\]
The fourth term must be estimated more accurately than in Section~\ref{sec:unicite}. Integrating by parts so as to make the derivative
of $w_1$ disappear, and then using $2ab\leq a^2+b^2$ on each of the two resulting terms,
\begin{align*}
    -2(J_1-J_2)\int_\R (w_1-w_2)\,\partial_x(u_\infty w_1)\varphi_X\,dx
    &= 2(J_1-J_2)\int_\R \big(\partial_x(w_1-w_2)\varphi_X+(w_1-w_2)\varphi_X'\big)u_\infty w_1 \,dx \\
    \leq{} &\int_\R\big(\partial_x(w_1-w_2)\big)^2u_\infty\varphi_X\,dx \\
    & +\big(J_1-J_2\big)^2\int_\R\frac{u_1^2}{u_\infty}
    \Big(\varphi_X+\frac{(\varphi_X')^2}{\varphi_X}\Big)dx
    +\underbrace{\int_\R (w_1-w_2)^2u_\infty\varphi_X\,dx}_{\leq\, I(t)} .
\end{align*}
The first term is absorbed by $-2\int_\R\big(\partial_x(w_1-w_2)\big)^2u_\infty\varphi_X\,dx$.
For the second one we use \eqref{est:borne L^infini J}, which gives
$u_1(t,\cdot)\leq \big \|\tfrac{u^0}{Q}\big\|_\infty e^{\Lambda_{J_1}(t)}\,Q$ on $(0,T)$, together
with $\Lambda_{J_1}(t)\leq C\big(T+R\sqrt{T}\big)$ for every $t\in(0,T)$, by the  Cauchy--Schwarz inequality.
Setting
\[
M(T,R) := \Big\|\frac{u^0}{Q}\Big\|_\infty^2\exp\Big(2C\big(T+R\sqrt{T}\big)\Big),
\]
we obtain
\[
    \int_\R\frac{u_1^2}{u_\infty}\Big(\varphi_X+\frac{(\varphi_X')^2}{\varphi_X}\Big)dx
    \;\leq\;M(T,R)\,\eta(X),
    \qquad
    \eta(X):=\int_\R\frac{Q^2}{u_\infty}\Big(\varphi_X+\frac{(\varphi_X')^2}{\varphi_X}\Big)dx .
\]
Notice that $\eta(X)\leq \Phi_0\int_X^{\infty}\frac{Q^2}{u_\infty}\,dx\to0$ as $X\to\infty$.
Indeed $Q\sim c_0/h$ by Lemma~\ref{lem:proprietes Q} and $u_\infty\sim N_\infty/h$ at $+\infty$,
so that $Q^2/u_\infty\sim c_0^2/(N_\infty h)$ is integrable near $+\infty$ by
\eqref{eq:behavior h +infty}.

Finally, integrating by parts in the remaining term,
\begin{align*}
    -2J_2\int_\R (w_1-w_2)\,\partial_x\big(u_\infty (w_1-w_2)\big)\varphi_X \,dx
    &=-J_2\int_\R (w_1-w_2)^2\big(u_\infty'\varphi_X-u_\infty \varphi'_X\big)\,dx\\
    &\leq |J_2(t)|\,\big(\kappa_\infty+\|\varphi'\|_\infty\big)\int_\R (w_1-w_2)^2u_\infty\,dx\\
    &= 2\,\theta_3\,|J_2(t)|\,I(t),
\end{align*}
where we used
$|u_\infty'|\leq\kappa_\infty u_\infty$ on $\operatorname{supp}\varphi_X\subset[X,\infty)$.

Gathering all the previous estimates and discarding the remaining non-positive diffusion term,
we obtain, for any $t\in(0,T)$,
\[
    I_X'(t)+N_\infty\Big(\frac{N_1(t)-N_2(t)}{N_\infty}\Big)^2
    \leq \Theta_X(t)\,I(t)
    +M(T,R)\,\eta(X)\,\big(J_1(t)-J_2(t)\big)^2,
\]
with $\Theta_X(t):=1+N_\infty\theta_1(X)+\theta_2(X)+2\theta_3(X)|J_2(t)|$. Integrating over
$(0,T)$ and using $I_X(0)=0$ and $I_X(T)\geq0$,
\begin{align}\label{eq:ineg N_1-N_2}
    \int_{0}^{T} N_{\infty}\Big(\frac{N_{1}-N_{2}}{N_{\infty}}\Big)^{2} dt
    \le {} & \int_{0}^{T}\Theta_X(t)\, I(t)\,dt
    +M(T,R)\,\eta(X)\int_{0}^{T}\big(J_{1}-J_{2}\big)^2\,dt .
\end{align}
To bound the first term we use \eqref{eq: controle I}, which gives, for any $t\in(0,T)$,
\begin{equation}\label{eq:borne I par R}
    I(t)\;\leq\;\|u^0\|_{L^2_{\mathrm w}}^{2}\,e^{2R^2}
    \int_{0}^{T}\big(J_{1}(\tau)-J_{2}(\tau)\big)^{2}\,d\tau ,
\end{equation}
together with $T\leq\sqrt T$ and $\|J_2\|_{L^2(0,T)}\leq R$,
\[
    \int_{0}^{T}\Theta_X(t)\,dt
    \;\leq\;\big(1+N_\infty\theta_1(X)+\theta_2(X)\big)T+2\theta_3(X)R\sqrt T
    \;\leq\;\Theta(X)\,(1+R)\sqrt{T},
\]
where $\Theta(X):=\max\big(1+N_\infty\theta_1(X)+\theta_2(X),\,2\theta_3(X)\big)$. Combining
these estimates yields
\[
\int_{0}^{T} N_{\infty}\Big(\frac{N_{1}-N_{2}}{N_{\infty}}\Big)^{2} dt
\leq \Big(\Theta(X)(1+R)\sqrt{T}\,\|u^0\|_{L^2_{\mathrm w}}^{2}\,e^{2R^2}
+M(T,R)\,\eta(X)\Big)\int_{0}^{T}\big(J_{1}-J_{2}\big)^2\,dt .
\]
The left-hand side does not depend on $X$, so we may take the infimum over $X>0$ and set
\[
  C(R,T,u^0):=\inf_{X>0}\Big(\Theta(X)(1+R)\sqrt{T}\,\|u^0\|_{L^2_{\mathrm w}}^{2}\,e^{2R^2}
  +M(T,R)\,\eta(X)\Big),
\]
which gives \eqref{eq:N-J lipschitz}. Finally, let $\varepsilon>0$. Since $M(T,R)$ is bounded as
$T\to0$ and $\eta(X)\to0$ as $X\to\infty$, we may fix $X_\varepsilon$ such that
$\sup_{T<1}M(T,R)\,\eta(X_\varepsilon)<\tfrac{\varepsilon}{2}$; then, $X_\varepsilon$ being
fixed, the first term is smaller than $\tfrac{\varepsilon}{2}$ for $T$ small enough. Hence
$\limsup_{T\to0}C(R,T,u^0)\leq\varepsilon$ for every $\varepsilon>0$, that is
$\lim_{T\to0}C(R,T,u^0)=0$.
\end{proof}

Equipped with this result, the global existence for an initial datum in $\mathcal{X}$ follows
from the Picard fixed point theorem.

\begin{prop}\label{prop:existence unicite glob u^0 dans X}
Let $u^0\in \mathcal{X}$ such that $u^0\geq 0$. Then there exists a unique global solution $(u,N)$ of
\eqref{eq:FPSL Non-lin} in $\mathcal{X}$.
\end{prop}

\begin{proof}
Uniqueness is an immediate consequence of Proposition~\ref{prop:unicite L2w}, so we only have to
prove existence.

We first construct a local solution. Given $J\in L^2_{\mathrm{loc}}(0,\infty)$, we denote by $(u,N)$ the solution of the linear
problem \eqref{eq:FPSL with J(t)} with initial datum $u^0$ provided by
Theorem~\ref{th:existence unicite avec J(t)}, and we set $\Phi(J):=b\,N$. By construction, a
pair $(u,N)$ solves \eqref{eq:FPSL Non-lin} on $(0,T)$ if and only if $bN$ is a fixed point of
$\Phi$. For $T,R>0$ we let $X_{T,R}:=\{f\in L^2(0,T):\|f\|_{L^2(0,T)}\le R\}$, so that
$J\in X_{T,R}$ means $\mathcal{J}(T)\leq R^2$.

Let $T\in(0,1)$ and $J\in X_{T,R}$. Then $\int_0^{T}|J|\leq\sqrt T\,\|J\|_{L^2(0,T)}\leq R$ by
Cauchy--Schwarz, so that \eqref{est:borne L^infini J}, together with $hQ\to c_0$ at $+\infty$,
gives the pointwise bound $N(t)\leq c_0\|u^0/Q\|_\infty e^{C(1+R)}$ on $(0,T)$, whence
\[
  \|\Phi(J)\|_{L^2(0,T)}^{2}=b^{2}\int_0^{T}N(t)^{2}\,dt
  \;\leq\;c_0^2\,b^{2}\Big\|\frac{u^0}{Q}\Big\|_\infty^{2}e^{2C(1+R)}\,T ,
\]
while, for $J_1,J_2\in X_{T,R}$, estimate \eqref{eq:N-J lipschitz} gives
\[
  \big\|\Phi(J_1)-\Phi(J_2)\big\|_{L^2(0,T)}^{2}
  \;\leq\;b^{2}C(R,T,u^0)\;\|J_1-J_2\|_{L^2(0,T)}^{2}.
\]
Fix any $R>0$ and choose $T_0\in(0,1)$ small enough so that
\begin{equation}\label{new:choix}
  c_0^2\,b^{2}\Big\|\frac{u^0}{Q}\Big\|_\infty^{2}e^{2C(1+R)}\,T_0\;\leq\;R^{2}
  \qquad\text{and}\qquad
  b^{2}C(R,T_0,u^0)<1 ,
\end{equation}
which is possible since $\lim_{T\to0}C(R,T,u^0)=0$. The first condition makes $X_{T_0,R}$
invariant, the second makes $\Phi$ contracting,
\begin{equation}\label{state:existence locale}
  \Phi\big(X_{T_0,R}\big)\subset X_{T_0,R}
  \qquad\text{and}\qquad
  \Phi\big|_{X_{T_0,R}}\ \text{is a contraction.}
\end{equation}
By Picard's fixed point theorem there is a unique $\overline J\in X_{T_0,R}$ with
$\Phi(\overline J)=\overline J$; the solution of \eqref{eq:FPSL with J(t)} associated with
$\overline J$ then solves \eqref{eq:FPSL Non-lin} on $(0,T_0)$.

\medskip
Let $\tau$ be the maximal existence time of the solution so obtained, and assume $\tau<\infty$. Applying
\eqref{est:borne L^1 J} with $J=bN$ and $\varepsilon=\frac{1}{2|b|}$ gives
\[
\int_0^{\tau} |N(t)|\,dt\;\leq\;2\big(\gamma_\varepsilon\tau+1\big)<\infty ,
\]
and \eqref{est:borne L^infini J}, again with $J=bN$, then yields
\[
\Big\|\frac{u(\tau)}{Q}\Big\|_\infty\;\leq\;
\Big\|\frac{u^0}{Q}\Big\|_\infty\exp\Big(C\tau+C|b|\int_0^{\tau}|N(t)|\,dt\Big)<\infty .
\]
Hence $u(\tau)\in\mathcal{X}$, and restarting from $u(\tau)$ extends the solution
strictly beyond $\tau$, a contradiction. Therefore $\tau=\infty$.
\end{proof}

\section{Global estimates for the nonlinear problem}
\label{sec:estimates}

We start by studying the stationary states associated with \eqref{eq:FPSL Non-lin} which are useful for several bounds of this section.
Then, we use $L^2-$entropy estimates and $L^2-$truncated entropy to prove Theorem~\ref{th:cesaro bound} and Theorem~\ref{th:cv exp}.
Finally, we prove Theorem~\ref{th:existence and uniqueness FPSLNL} in its full generality using the above estimates.

\subsection{Stationary states}
\label{sec:stationary states}

Let $b\in \R$. We study the properties of the normalized stationary states $(u_\infty^b,N_\infty^b)$. Depending on the value of the minimum of the drift function $h$, one or several stationary states may exist. Hence, in order to distinguish these cases, we sometimes use the form
\begin{equation} \label{h=h0+A}
h(x) =h_0(x)-A \quad \text{with}\quad  A\in\R, \quad \min_\R h_0 = 0.
\end{equation}

\begin{prop}[Stationary states and connectivity $b$]\label{prop:stationary states}
For all $b\in \R$, there is at least one normalized stationary state of \eqref{eq:FPSL Non-lin} and
\\
1. it is unique away from a possibly empty interval $[b_-, b_+]$ with $0<b_- <b_+< \infty$,  and for $b<0$, $N_\infty^b$ remains bounded,  
\\
2. when $\min_\R h \geq 0$, this interval is empty,
\\
3. with the form \eqref{h=h0+A} and $A$ large enough, for some $b$'s  there are at least three stationary states,
\\
4. when $h$ satisfies also \eqref{eq:croissance h Poincaré}, then for  $b\in \R$, the stationary state(s) $u_\infty^b$ are weights that allow the Poincaré inequality \eqref{eq:Poincaré uinfty} with a coefficient $\alpha_b>0$ which is bounded on any interval $(-B, B]$ with $B>0$.
\end{prop}

\begin{figure}[htbp]
\centering
\includegraphics[width=.9\textwidth]{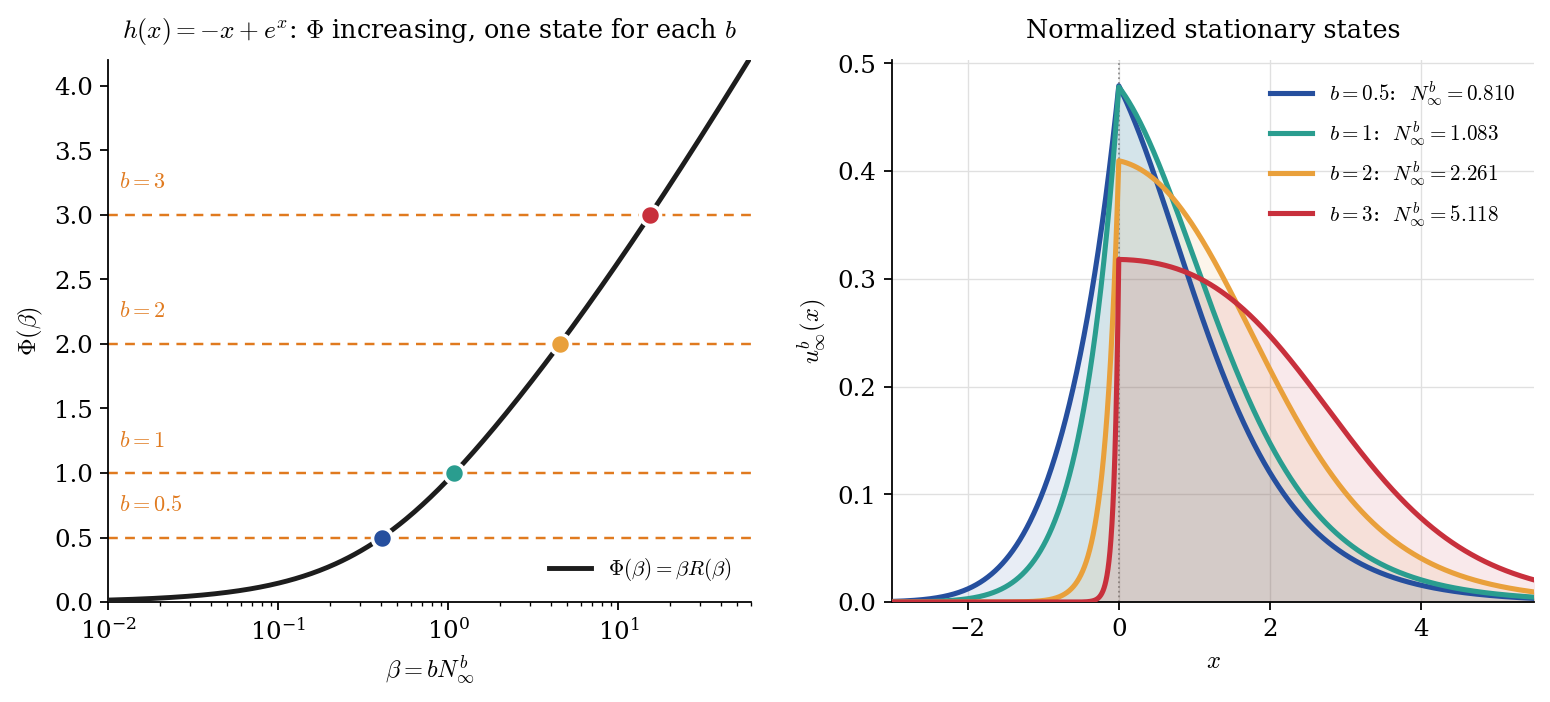}
\vspace{-15pt}
\caption{Drift $h(x)=-x+e^{x}$. Left: $\Phi$ is increasing, so there is exactly one stationary state for each $b$. Right: the normalized stationary states for $b\in\{0.5,1,2,3\}$.}
\label{fig:unique state}
\end{figure}

\begin{figure}[htbp]
\centering
\includegraphics[width=.9\textwidth]{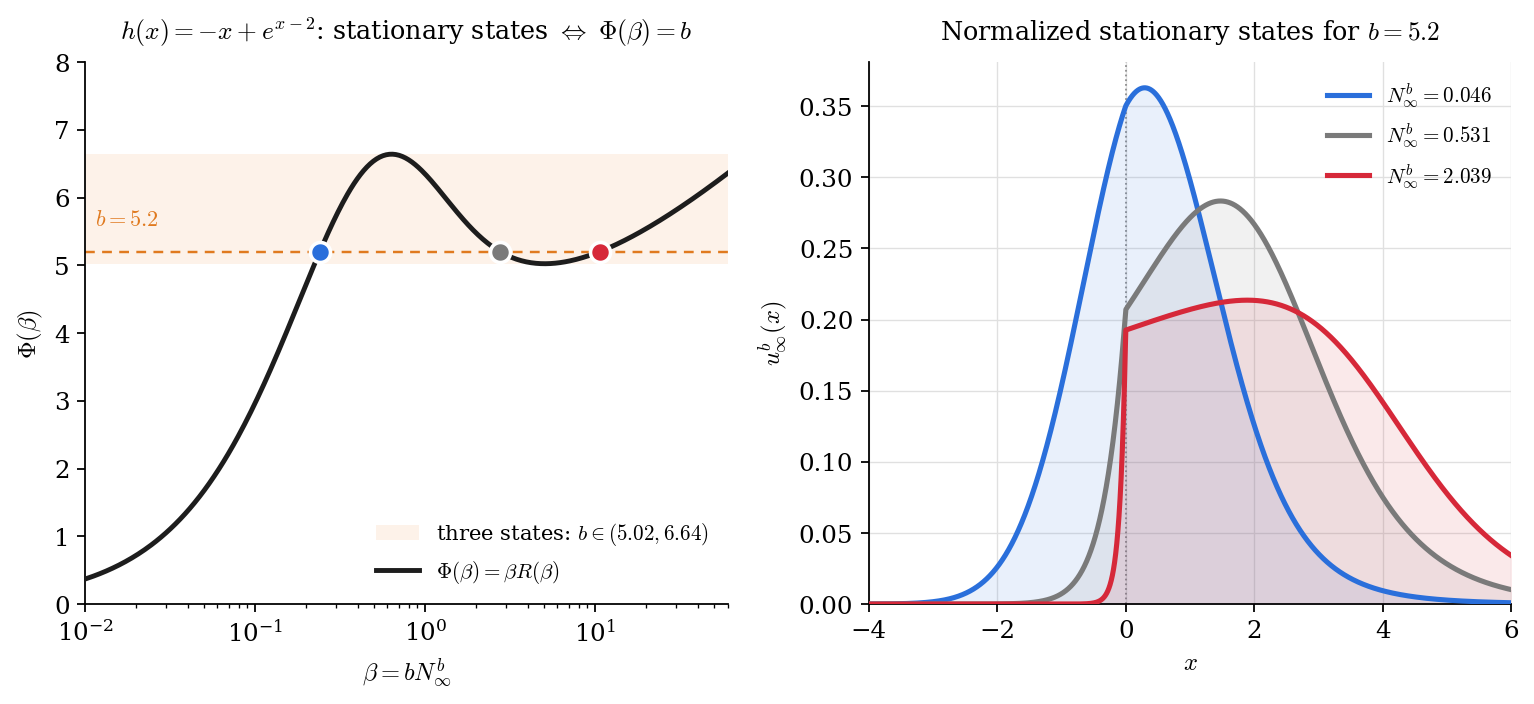}
\vspace{-15pt}
\caption{Drift $h(x)=-x+e^{x-2}$. Left: the map $\Phi(\beta)=\beta R(\beta)$; stationary states correspond to the solutions of $\Phi(\beta)=b$, and there are three of them for $b\in(5.02,6.64)$. Right: the three normalized stationary states for $b=5.2$.}
\label{fig:three states}
\end{figure}

These results are illustrated in Fig.~\ref{fig:unique state} and~\ref{fig:three states}.

\begin{proof}
After integration of the equation for $u_\infty^b$, we find
\begin{equation}\label{eq:equation u'_b}
     \big(h(x)+bN_\infty^b\big) u_\infty^b(x)-(u_\infty^b)'(x)=N_\infty^b\mathbf{1}_{x\geq 0},
     \quad N_\infty^b := \lim_{x \rightarrow \infty} h(x)u^b_\infty(x),
     \quad \int_{-\infty}^{+\infty} u^b_\infty(x) dx=1 .
\end{equation}
Solving this differential equation, we get for all \(x \in \R, \)
\begin{equation}\label{eq:ecriture u_b}
 u_\infty^b(x)=N_\infty^b e^{H(x)+bN_\infty^b x} 
 \int_{\max(0,x)}^\infty e^{-H(y)-bN_\infty^b y }\,dy, \qquad H(x):=\int_0^x h(y)dy.
\end{equation}
Notice that in this expression the connectivity $b$ and the activity $N_\infty^b$ only appear through their product. We may therefore define, for $\beta\in\R$,
\begin{equation}\label{def:R}
    R(\beta) :=\int_\R\exp\Big(H(x)+\beta x\Big)
    \int_{\max(0,x)}^\infty\exp\Big(-H(y)-\beta y\Big)\,dy\,dx,
\end{equation}
which is the inverse activity of the linear stationary problem with frozen drift $h+\beta$.
The mass constraint in \eqref{eq:equation u'_b} then leads to the condition \( 1= N_\infty^b R\big(bN_\infty^b\big) \) which means we need to find $\beta$ such that
\begin{equation}\label{eq:point fixe}
    b = \beta R(\beta):=\Phi(\beta), \qquad \beta = b N_\infty^b.
\end{equation}

Exchanging the order of integration in \eqref{def:R} and setting $z=y-x$, we obtain
\begin{equation}\label{eq:R Laplace}
\begin{cases}
R(\beta)
%=\int_0^\infty\!\!\int_{-\infty}^y e^{H(x)-H(y)-\beta(y-x)}\,dx\,dy
=\int_0^\infty e^{-\beta z}K(z)\,dz,
\qquad K(z):=\int_0^\infty e^{H(y-z)-H(y)}\,dy >0,
\\
R'(\beta)<0, \qquad R''(\beta)>0,
\end{cases}
\end{equation}
because $K$ is positive and $C^1$ on $(0,\infty)$. Furthermore, we have 
\[
K'(z)=-\int_0^\infty h(y-z)\,e^{H(y-z)-H(y)}\,dy ,
\]
and $K$ satisfies
\begin{equation}\label{eq:proprietes K}
\lim_{z\to0^+}K(z)=+\infty,\qquad \lim_{z\to0^+}K'(z)=-\infty,\qquad 
\exists\gamma,z_1>0,\ \ K(z)\leq C_\gamma e^{-\gamma z}\ \text{ for } z\geq z_1 ,
\end{equation}
The first 2 points come from the inequalities:
\[
K(0)\geq  \sup_{C\geq 0}\lim_{z\rightarrow0}\int_0^C e^{H(y-z)-H(y)}\,dy = +\infty \quad \text{and}\quad -K'(0)\geq  \sup_{C\geq 0}\lim_{z\rightarrow 0}\int_0^C h(y-z)e^{H(y-z)-H(y)}\,dy = +\infty.
\]
In addition, because of the linear tail of $h$ at $-\infty$ stated in \eqref{eq:behavior h -infty}, we can find $x_2<0$ such that $h(x)\geq \frac{1}{2}x$ for $x\leq x_2$ and the superlinearity of $h$ on $+\infty$, we can find $\gamma >0$ and $z_1>0$ big enough such that  we have
\[
\fa y\in \R, \quad \fa z>z_1 \qquad\int_{y-z}^yh(w)\,dw \geq \gamma z(y+1).
\]
This ensures that $\fa z>z_1$, we have
\[
K(z)\leq \frac{e^{-\gamma z}}{\gamma z_1}.
\]
Now, because of the properties of $K$ near $0$, we can find $z_0>0$ such that  
\begin{equation}\label{eq:croissance K}
    K \text{ is decreasing on }(0,z_0) \quad \text{and}\quad K(z_0)=\sup_{z\geq z_0}K(z).
\end{equation}

These properties \eqref{eq:proprietes K} imply that there are $0<\beta_-<\beta_+$ such that $\Phi$ satisfies
\begin{equation}\label{eq:proprietes Phi}
\Phi'(\beta) >0, \quad \Phi''(\beta) <0 \quad \text{for}\quad \beta \leq \beta_-, \qquad \Phi'(\beta) >0, \quad \text{for}\quad \beta \geq \beta_+.
\end{equation}

Using the first point of \eqref{eq:proprietes K}, we deduce
\[
\lim_{\beta\rightarrow\infty}\Phi(\beta) = \lim_{\beta\rightarrow\infty} \beta\int_0^\infty e^{-\beta z}K(z)\,dz = +\infty.
\]
Therefore, by writing $b_-=\max \Phi^{-1}((-\infty,\beta_-])$ and $b_+=\min \Phi^{-1}([\beta_+,\infty))$, the equation $b=\Phi(bN)$ admits a unique solution $N_\infty^b$ for every $b\in(-\infty,b_-)\cup (b_+,\infty)$. Using the concavity of $\Phi$ on $(-\infty,0]$, we obtain that for any $\beta\leq 0$,
\[
\Phi(\beta)\leq \beta R(0) \quad \text{and}\quad
b =\Phi(bN_\infty^b)\leq bN_\infty^b R(0).
\]
hence, we get 
\[
N_\infty^b\leq \frac{1}{R(0)}.
\]
The first item of Proposition~\ref{prop:stationary states} follows immediately.

To prove \eqref{eq:proprietes Phi}, we first notice that 
\[
\Phi'(\beta)=\beta R'(\beta)+R(\beta) = \int_0^\infty(1-\beta z)e^{-\beta z}K(z)\,dz,
\]
and using the fact that $\int_0^\infty (1-\beta z)e^{-\beta z}\,dz=0$, we obtain that for any $\beta >0$, 
\[
\Phi'(\beta)= \int_0^\infty(1-\beta z)e^{-\beta z}\big(K(z)-K(\frac{1}{\beta})\big)\,dz,
\]
Therefore, using \eqref{eq:croissance K}, we have for $\beta \geq \frac{1}{z_0}$ 
\[
\fa z>0, \quad (1-\beta z) \big(K(z)-K(\frac{1}{\beta})\big)\geq 0, \qquad \Phi'(\beta)>0.
\]
By taking $\Phi(b^+)= \max_{0 \leq \beta \leq \frac{1}{z_0}} \Phi(\beta)$, this ensures that for $b>b^+$, the solution $b=\Phi(\beta)$ has only one solution and the first point is proved.

For the second point, if we assume in addition that $\min_\R h\geq 0$, then $K'(w)<0$ for any $w>0$. By the above argument, $\Phi$ is strictly increasing on $(0,\infty)$ and therefore, the stationary state is unique for any $b\in \R$.
\\

Next we turn to point 3. using the representation \eqref{h=h0+A}. Then, it follows that $H(x) = -Ax+H_0(x)$ with $H_0(x)=\int_0^x h_0$ and $K(z):=\int_0^\infty e^{Az+H_0(y-z)-H_0(y)}\,dy=e^{Az}K_0(z)$. Setting $R_0(\beta) = \int_0^\infty e^{-\beta z}K_0(z)\,dz$ we can write 
\[
\Phi(\beta) = \beta \int_0^\infty e^{-(\beta-A)z}K_0(z)\,dz=  \beta R_0(\beta-A) .
\]
Because $\Phi$ is positive, increasing for $\beta >  0$ small or large enough, it is sufficient to find  $\beta_A>0$ such that $\Phi'(\beta_A)<0$ in order to get non-uniqueness of steady states. To do so, we observe that
\[
 \Phi'(\beta) = R_0(\beta-A)+\beta R_0'(\beta-A), \qquad \Phi'(A) = AR_0'(0)+R_0(0).
\]
Because, see \eqref{eq:R Laplace}, $R_0'(0)<0$, for $A$ large enough $\Phi'(A)<0$ and the point 3. is proved.
\\

 Finally, in the case where $h$ satisfies the growth condition
\eqref{eq:croissance h Poincaré}, for any $b\in \R$ the function
$x\mapsto h(x)+bN_\infty^b$ satisfies also \eqref{eq:croissance h Poincaré}. Then we are
reduced to the case $b=0$, and thus $\alpha_b<\infty$.  By the Muckenhoupt criterion in dimension one \cite{muckenhouptHardysInequalityWeights1972}, $\alpha_b\leq 4A(b)$, where
\[
  A(b)=\max\bigg(
  \sup_{x>0}\Big(\int_x^{\infty}u^b_\infty\Big)\Big(\int_{0}^{x}\frac{dy}{u^b_\infty}\Big),\;
  \sup_{x<0}\Big(\int_{-\infty}^{x}u^b_\infty\Big)\Big(\int_{x}^{0}\frac{dy}{u^b_\infty}\Big)
  \bigg).
\]
One easily checks that $b\mapsto A(b)$ is locally bounded.

\end{proof}

The function $\Phi$ is plotted in Figures~\ref{fig:unique state} and~\ref{fig:three states} for two examples of drift, $h(x)=-x+e^{x}$ and $h(x)=-x+e^{x-2}$. In the first case, the associated function $\Phi$ is strictly increasing, as established in point~2 of Proposition~\ref{prop:stationary states}, and Figure~\ref{fig:unique state} shows the density of the stationary state for several values of $b>0$. The second case, shown in Figure~\ref{fig:three states}, illustrates point~3: there is an interval of values of $b$ for which three stationary states exist. The numerical results also suggest that the middle stationary state is unstable, while the two others are stable.

More generally, the deterministic dynamics suggests that the number of stationary states for a given function $h$ is related to the number of zeros of $h$.

\subsection{$L^2$-Cesàro bounds on the activity $N(t)$ for inhibitory connections and exponential convergence in the weakly nonlinear regime}

We now establish global $L^2_{\mathrm w}(\R)$ bounds on the solution of \eqref{eq:FPSL Non-lin}. With this established, we conclude the exponential convergence towards the stationary state in 
the case where $|b|$ is small enough and $u_\infty$ satisfies the Poincar\'e inequality~\eqref{eq:Poincaré uinfty}. Since the sign of $b$ plays 
a key role in the $L^2$-Cesàro bounds of $N(t)$, we split the proof of Theorem \ref{th:cv exp} in two parts.

 Our proofs make heavy use of the entropy equality associated with  the nonlinear equation relative to a stationary state $u_\infty^{b'}$ associated with a connectivity $b'\in\R$ to be chosen. For a solution of \eqref{eq:FPSL Non-lin}, let $w_{b'}(t) = \frac{u(t)}{u_\infty^{b'}}$, we write two equalities that we use extensively later.
\\
1. \textbf{Entropy equality with the nonlinear term, with respect to the stationary state $u_\infty^{b'}$}.
\begin{align}\label{eq:entropienonlin}
\frac{d}{dt}\int_\R H\big(w_{b'}(t)\big) u_\infty^{b'}\,dx 
&= -\int_\R H''\big(w_{b'}(t)\big) \big(\partial_x w_{b'}(t)\big)^2  u_\infty^{b'}\,dx \notag
\\
& -N_\infty^{b'}\bigg(H\Big(\frac{N(t)}{N_\infty^{b'}}\Big)-H\big(w_{b'}(t,0)\big)-H'\big(w_{b'}(t,0)\big)\Big(\frac{N(t)}{N_\infty^{b'}}-w_{b'}(t,0)\Big)\bigg) \notag\\
&-\big(bN(t)-b'N_\infty^{b'}\big)\int_\R \Big(w_{b'}(t)H'(w_{b'}(t))-H(w_{b'}(t))\Big)\big(u_\infty^{b'}\big)'\,dx.
\end{align}

2. \textbf{Truncated entropy equality with the nonlinear term, with respect to the stationary state $u_\infty^{b'}$}.
Recall the cut-off function  $\varphi$ defined on~\eqref{def:varphi}.
\begin{align}\label{eq:entropietruncnonlin}
\frac{d}{dt}\int_{\R} H(w_{b'}(t))\,u_\infty^{b'}\varphi\,dx
&= -\int_{\R} \big(\partial_x w_{b'}(t)\big)^2 H''(w_{b'}(t))\,u_\infty^{b'}\varphi\,dx 
%\notag\\ &\quad 
- N_\infty^{b'}\!\left(
    H\!\left(\frac{N(t)}{N_\infty^{b'}}\right)
    - \int_{\R} H(w_{b'}(t))\,\varphi'\,dx
\right) \notag\\
&\quad - \big(b N(t) - b' N_\infty^{b'}\big)
\int_{\R} \Big[
    w_{b'}(t)\,H'(w_{b'}(t))\,(u_\infty^{b'})'\,\varphi
    - H(w_{b'}(t))\,(u_\infty^{b'}\varphi)'
\Big]\,dx \notag\\
&\quad + \int_{\R} H(w_{b'}(t))\,\big(u_\infty^{b'}\varphi'\big)'\,dx.
\end{align}

\begin{prop}[Uniform $L^2$-Cesàro bound in the inhibitory case]
\label{prop:cesaro inhib}
Let $b\leq0$. There is a constant $C_\infty>0$,
depending only on $u_\infty$, such that the
solution $(u,N)$ of \eqref{eq:FPSL Non-lin} with a probability density initial datum $u^0\in L^2_{\mathrm w}(\R)$, satisfies,  for all $T\geq 0$
\begin{equation}\label{eq: cesar L^2 inhib}
     \int_{\R} w(T)^2\,u_\infty\,\varphi\,dx
    +\int_0^T \frac{N(t)^2}{N_\infty}\,dt\leq \int_{\R}\!\left(\frac{u^0}{u_\infty}\right)^{\!2}u_\infty\,\varphi\,dx + C_\infty T.
\end{equation}
\end{prop}

\begin{proof}
We start from equation \eqref{eq:entropietruncnonlin} with $H(x)=x^2$ and $b'=0$. Consistently with the previous notation, we write $u_\infty=u_\infty^0$ and
$N_\infty=N_\infty^0$ and get  
\begin{align*}
\frac{d}{dt}\int_{\R} w(t)^2\,u_\infty\varphi\,dx
&= -2\int_{\R} \big(\partial_x w(t)\big)^2\,u_\infty\varphi\,dx - N_\infty\bigg(
    \Big(\frac{N(t)}{N_\infty}\Big)^{\!2}
    - \underbrace{\int_{\R} w(t)^2\,\varphi'\,dx}_{=:\,I_1}
\bigg) \\
&\quad - b N(t)
\int_{\R} \Big[
    2w(t)^2\,(u_\infty)'\,\varphi
    - w(t)^2\,(u_\infty\varphi)'
\Big]\,dx
+ \underbrace{\int_{\R} w(t)^2\,\big(u_\infty\varphi'\big)'\,dx}_{=:\,I_2}.
\end{align*}
Notice that, since $b\leq 0$, $\varphi'\geq 0$ and $u_\infty'\leq0$  on $\operatorname{supp}\varphi\subset[x_1,\infty)$,
we have
\begin{align*}
- b N(t)\int_{\R} \Big[2w(t)^2\,(u_\infty)'\,\varphi
    - w(t)^2\,(u_\infty\varphi)'\Big]\,dx
&= -b N(t)\int_{\R} w(t)^2\,\big(u_\infty'\,\varphi - u_\infty\,\varphi'\big)\,dx \leq 0.
\end{align*}

We now estimate $I_1$ and $I_2$. Let $K=\operatorname{supp}(\varphi')$ and set
$C_K=\left\|\frac{1}{u_\infty}\right\|_{\infty,K}$. We have
\begin{align}
I_1 = \int_{\R} w(t)^2\,\varphi'\,dx
&\leq C_K\int_{\R} w(t)^2\,u_\infty\,\varphi'\,dx
= C_K\int_{\R} w(t)\,u(t,x) \,\varphi'\,dx\notag
\\
&= C_K\int_{\R} w(t,x)\,\partial_x\!\left(\int_{-\infty}^x\!u(t,y)\,dy\right)\varphi'(x)\,dx \notag\\
&= -C_K\int_{\R} \left(\int_{-\infty}^x\!u(t,y)\,dy\right)
    \!\big(\partial_x w(t,x)\,\varphi'(x)+w(t,x)\,\varphi''(x)\big)\,dx \notag\\
&\leq C_K\int_{\R}\Big(|\partial_x w(t,x)|\,\varphi'(x)+w(t,x)\,|\varphi''(x)|\Big)\,dx,
\end{align}
where in the last step we used mass conservation \eqref{eq:conservation de la masse}, 
$\int_{-\infty}^{x}u(t,y)\,dy \leq 1$. Since
$\operatorname{supp}(\varphi'')\subset K$, we have
\[
\int_{\R} w(t)\,|\varphi''|\,dx
= \int_{\R} u(t)\,\frac{|\varphi''|}{u_\infty}\,dx
\leq \left\|\frac{\varphi''}{u_\infty}\right\|_{\infty}.
\]
Furthermore, for any $\varepsilon>0$, using that $|ab|\leq \varepsilon a^2+\frac{1}{4\varepsilon}b^2$, we obtain
\[
\int_{\R} |\partial_x w(t)|\,\varphi'\,dx
\leq \varepsilon\int_{\R} \big(\partial_x w(t)\big)^2 u_\infty\,\varphi\,dx
+ \frac{1}{4\varepsilon}\left\|\frac{(\varphi')^2}{u_\infty\varphi}\right\|_{1}.
\]
Let $\psi=\frac{(\varphi')^2}{\varphi}$. Notice that $\psi\in\mathcal{C}_c^\infty(\R)$, we obtain the following bound on $I_1$,
\begin{equation}\label{eq:estimI1}
I_1 \leq C_K\!\left(
    \varepsilon\int_{\R}\big(\partial_x w(t)\big)^2 u_\infty\,\varphi\,dx
    + \frac{1}{4\varepsilon}\left\|\frac{\psi}{u_\infty}\right\|_{1}
    + \left\|\frac{\varphi''}{u_\infty}\right\|_{\infty}
\right).
\end{equation}
In the same spirit,
\begin{align*}
    I_2 =& \int_\R w(t)^2\big(u_\infty'\varphi'+u_\infty\varphi'')\,dx\\
    &\leq\big( \|u_\infty'\|_\infty+ \|u_\infty\|_\infty\big)\int_\R w(t)^2(\varphi'+|\varphi''|)\,dx
\end{align*}
and by replacing $\varphi'$ by $|\varphi''|\in \mathcal{C}^\infty_c(\R)$ in \eqref{eq:estimI1}, we obtain the existence of $C_K, C(\varepsilon)>0$ such that 
\[
N_\infty I_1+I_2\leq C_K\varepsilon\int_\R \big(\partial_xw(t)\big)^2u_\infty \varphi\,dx +C(\varepsilon)
\]

we obtain
\[
\frac{d}{dt}\int_{\R} w(t)^2\,u_\infty\,\varphi\,dx
+ \frac{N(t)^2}{N_\infty}
\leq -\left(2 -  C_K\varepsilon \right)
    \int_{\R}\big(\partial_x w(t)\big)^2 u_\infty\,\varphi\,dx + C(\varepsilon).
\]
Choosing $\varepsilon=\frac{1}{2 C_K}$ and $C_\infty:=C(\frac{1}{2 C_K})$, so that the bracket equals $\frac 3 2>0$,
and integrating over $(0,T)$ for any~$T\geq 0$, we conclude that
\begin{equation}\label{eq: ces L^2 inhib}
    \int_\R w(T)^2 u_\infty\,\varphi\,dx + \int_0^T \frac{N(t)^2}{N_\infty}\,dt
    \leq \int_{\R}\!\left(\frac{u^0}{u_\infty}\right)^{\!2}u_\infty\,\varphi\,dx + C_\infty T.
\end{equation}
This proves \eqref{eq: cesar L^2 inhib}.
\end{proof}

We now turn to the exponential convergence of the solution towards the stationary state in the
weakly connected inhibitory regime. The $L^2$-Cesàro bound just established, whose constants are
universal, is what makes the threshold $b^*_i$ independent of the initial datum.

\begin{proof}[Proof of Theorem~\ref{th:cv exp}~\rm(i)]
According to \eqref{eq:entropienonlin}, with $H(x)=x^2$ and $b=b'$, the integral
\( \mathcal{I}_b(t) = \displaystyle\int_\mathbb{R}\!\left(\frac{u(t)}{u_\infty^b}-1\right)^{\!2}u_\infty^b\,dx\)
satisfies
\[
\frac{d}{dt}\mathcal{I}_b(t)
= -2\int_\mathbb{R}(\partial_x w_{b}(t))^2\,u_\infty^b\,dx
  - N_\infty^b\left(\frac{N(t)}{N_\infty^b}-w_{b}(t,0)\right)^{\!2}
  - b\bigl(N(t)-N_\infty^b\bigr)\int_\mathbb{R} w_{b}(t)^2\,(u_\infty^b)'\,dx.
\]
To treat the last term, we integrate by parts and use that, because \( \int_\mathbb{R}w_{b}(t) u_\infty^b\,dx=1\),
\[
 \int_\mathbb{R}\big(w_{b}(t)\big)^2\,u_\infty^b\,dx =\int_\mathbb{R}\big((w_{b}(t)-1)^2+1\big)\,u_\infty^b\,dx =\mathcal{I}_b(t)+1.
\]
This implies 
\begin{align*}
b\bigl(N(t)-N_\infty^b\bigr)\int_\mathbb{R} w_{b}(t)^2\,(u_\infty^b)'\,dx
&= -2b\bigl(N(t)-N_\infty^b\bigr)
   \int_\mathbb{R}\partial_x w_{b}(t)\, w_{b} \,u_\infty^b\,dx 
   \\
&= \int_\mathbb{R}(\partial_x w_{b}(t))^2\,u_\infty^b\,dx
      + b^2\bigl(N(t)-N_\infty^b\bigr)^2
        \int_\mathbb{R}\big((w_{b}(t)-1)^2+1\big)\,u_\infty^b\,dx\\
&=\int_\mathbb{R}(\partial_x w_{b}(t))^2\,u_\infty^b\,dx+ b^2\bigl(N(t)-N_\infty^b\bigr)^2\big(\mathcal{I}_b(t)+1\big).
\end{align*}
 
Hence, we have
\[
\mathcal{I}_b'(t)
\leq -\int_\mathbb{R}(\partial_x w_{b}(t))^2\,u_\infty^b\,dx
     + (b\,N_\infty^b)^2\left(\frac{N(t)}{N_\infty^b}-1\right)^{\!2}\mathcal{I}_b(t)
     - N_\infty^b\left(\frac{N(t)}{N_\infty^b}-w_{b}(t,0)\right)^{\!2}+(bN_\infty^b)^2\left(\frac{N(t)}{N_\infty^b}-1\right)^{\!2}.
\]
For any $\varepsilon\in\bigl(0,\tfrac{N^b_\infty}{2}\bigr]$, applying $(x-y)^2\geq\theta(x-1)^2-2\theta(y-1)^2$
with $\theta=\varepsilon/N^b_\infty\leq\frac12$, we get
\[
-N^b_\infty\left(\frac{N(t)}{N_\infty^b}-w_{b}(t,0)\right)^{\!2}
\leq -\varepsilon\left(\frac{N(t)}{N_\infty^b}-1\right)^{\!2}+ 2\varepsilon\,(w_{b}(t,0)-1)^2 .
\]
Moreover, by the embedding of $H^1(-1,1)$ into $L^\infty(-1,1)$ (see \cite[Theorem~8.8]{brezisFunctionalAnalysisSobolev2011}), there exists $\overline{C}\geq 0$ such that
\[
(w_{b}(t,0)-1)^2 \leq \overline{C}\displaystyle\int_{-1}^{1}\Bigl(\bigl(\partial_x w_{b}(t,x)\bigr)^2+\bigl(w_b(t,x)-1\bigr)^2\Bigr)dx
\leq \frac{C_b'}{2}\Bigl(\displaystyle\int_\mathbb{R}\bigl(\partial_x w_{b}(t,x)\bigr)^2\,u_\infty^b\,dx+\mathcal{I}_b(t)\Bigr),
\]
where $C'_b:= 2 \overline{C}\bigl\|1/u_\infty^b\bigr\|_{\infty,(-1,1)}$. Hence
\[
\mathcal{I}_b'(t)
\leq -(1-C_b'\varepsilon)\int_\mathbb{R}(\partial_x w_{b}(t))^2\,u_\infty^b\,dx
     +\varepsilon C_b'\,\mathcal{I}_b(t)
     + \left(\frac{N(t)}{N_\infty^b}-1\right)^{\!2}
       \!\left((N_\infty^b\,b)^2\,\mathcal{I}_b(t) +(bN_\infty^b)^2 -\varepsilon\right).
\]

Using point 4 of Proposition~\ref{prop:stationary states}, the Poincaré inequality \eqref{eq:Poincaré uinfty} applied to $w_b(t)$ gives
\[
  \mathcal{I}_b(t)=\int_\R\big(w_b(t)-1\big)^2u^b_\infty\,dx
  \;\leq\;\alpha_b\int_\R\big(\partial_xw_b(t)\big)^2u^b_\infty\,dx.
\]
Therefore, we obtain
\begin{equation}\label{eq: I' joli}
\mathcal{I}_b'(t)\leq -\nu_b(\varepsilon)\,\mathcal{I}_b(t)
-\Bigl(\varepsilon-(N_\infty^b\,b)^2\bigl(1+\mathcal{I}_b(t)\bigr)\Bigr)\Big(\frac{N(t)}{N_\infty^b}-1\Big)^{\!2},
\qquad
\nu_b(\varepsilon):=\frac{1-C'_b\varepsilon}{\alpha_b}-C'_b\varepsilon.
\end{equation}
We choose 
\begin{equation}\label{eq:minoration nu}
\varepsilon= \varepsilon_b:=\min\Big(\frac{N^b_\infty}{2},\frac{1}{4C'_b\bigl(1+\alpha_b\bigr)}\Big), \qquad  \nu_b(\varepsilon_b)\geq  \frac{1}{2\alpha_b} .
\end{equation}

 By Proposition~\ref{prop:stationary states}, we have $(bN^b_\infty)^2\to0$ and $\varepsilon_b\geq \varepsilon_0>0$ in the vicinity of $0$.
 Hence there is $b_i^*<0$, depending only on $h$, such that $(bN^b_\infty)^2\leq\frac{\varepsilon_0}{4}$
and $\varepsilon_b\geq\frac{3\varepsilon_0}{4}$ for every $b\in(b_i^*,0)$. With $\eta:=\frac{\varepsilon_0}{2}$,
\eqref{eq: I' joli} and \eqref{eq:minoration nu} give
\begin{equation}\label{eq: mathcal I}
\fa t\geq0,\qquad \mathcal{I}_b'(t)\leq -\Big(\frac{1}{2\alpha_b}-(N_\infty^b\,b)^2\Big(\frac{N(t)}{N_\infty^b}-1\Big)^{\!2}\Big)\mathcal{I}_b(t)-\eta\Big(\frac{N(t)}{N_\infty^b}-1\Big)^{\!2}.
\end{equation}

By Grönwall's lemma and the $L^2$-Cesàro bound \eqref{eq:cesaro L2}, which gives
$\int_0^t\big(N-N_\infty^b\big)^2\leq 2C_1(u^0)+A_b\,t$ with
$A_b:=2\big(C_2+(N_\infty^b)^2\big)$, we obtain
\begin{equation}\label{eq:gronwall inhib}
  \fa t \geq 0, \qquad \mathcal{I}_b(t)\;\leq\;
  \exp\Big(-\Big(\tfrac{1}{8\alpha_b}-b^2A_b\Big)\,t+2b^2C_1(u^0)\Big)\,\mathcal{I}_b(0).
\end{equation}

Again, we can assume $|b_i^*|>0$ small enough such that
\[
  \fa b\in(b_i^*,0),\qquad 2\lambda_b:=\frac{1}{8\alpha_b}-b^2A_b\;>\;0 .
\]
Define the time
\[
  T_b(u^0):=\frac{2b^2C_1(u^0)}{\lambda_b}\;\geq\;0 .
\]
For $t\geq T_b(u^0)$ one has $2b^2C_1(u^0)\leq \lambda_b t$, so that
\eqref{eq:gronwall inhib} gives
\[
  \fa t\geq T_b(u^0), \qquad \mathcal{I}_b(t)\;\leq\;e^{-\lambda_b t}\,\mathcal{I}_b(0),
\]
which proves the first inequality in \eqref{eq:cv exp N inhib}.

Finally, to show the second inequality in \eqref{eq:cv exp N inhib}, we consider the previous inequality \eqref{eq: mathcal I}

\[
(\eta-(N_\infty^b b)^2\mathcal{I}_b(t)) \left(\frac{N(t)}{N_\infty^b}-1\right)^{\!2}\leq -\mathcal{I}'_b(t)-\frac{1}{2 \alpha_b}\mathcal{I}_b(t)
\]
 
 Set $T^\star_b(u^0):=\max\Big(T_b(u^0),\ \frac{1}{\lambda_b}\log^+\frac{2(bN_\infty^b)^2\,\mathcal{I}_b(0)}{\eta}\Big)$.
For $t\geq T^\star_b(u^0)$, the first inequality in \eqref{eq:cv exp N inhib} gives
$(bN^b_\infty)^2\mathcal I_b(t)\leq\frac{\eta}{2}$, so \eqref{eq: mathcal I} and $\lambda_b\leq\frac{1}{2\alpha_b}$ yield
\[
\frac{\eta}{2}\Big(\frac{N(t)}{N_\infty^b}-1\Big)^{\!2}\leq-\big(\mathcal I_b'(t)+\lambda_b\,\mathcal I_b(t)\big).
\]
Multiplying by $e^{\lambda_bt}$ and integrating over $(T^\star_b,\infty)$, we obtain
\[
\frac{\eta}{2}\int_{T^\star_b}^{\infty}\Big(\frac{N(t)}{N_\infty^b}-1\Big)^{\!2}e^{\lambda_bt}\,dt
\leq\mathcal I_b(T^\star_b)\,e^{\lambda_bT^\star_b}.
\]

Since $N\in L^2_{\mathrm{loc}}(\R^+)$, the integral over $(0,T^\star_b)$ is finite as well, which gives the second inequality in \eqref{eq:cv exp N inhib}. This concludes the proof of Theorem~\ref{th:cv exp}~\rm(i).
\end{proof}

The excitatory case is different, because the argument used to establish the $L^2$-Cesàro bound
on $N(t)$ no longer applies. We therefore only establish exponential convergence, again in the
weakly connected regime, and the threshold on $b$ is now obtained by propagating in time the
smallness of $(N_\infty^b\,b)^2\bigl(1+\mathcal{I}_b(t)\bigr)$, which makes it depend on the initial
datum. We assume that $h$ satisfies \eqref{eq:croissance h Poincaré} and $u^0\in L^2_w\cap L^2_{w,b^+}(\R)$, then 
%for $b\in[0,b_-]$}
%\begin{equation}\label{eq:def gamma}\Gamma(b) := \frac{\bigl(N_\infty^b\,b\bigr)^{2}\bigl(1+\mathcal{I}_b(0)\bigr)}{\varepsilon_b},\qquad \mathcal{I}_b(0)=\int_\mathbb{R}\Bigl(\frac{u^0}{u_\infty^b}-1\Bigr)^{\!2}u_\infty^b\,dx,\end{equation}
%where $\varepsilon_b$ is the threshold of \eqref{eq:minoration nu}. Note that since we assume  then we ensure that $\mathcal{I}_b(0)$ is finite and $b\mapsto \mathcal{I}_{b}(0)$ is continuous, so is $\Gamma$ on $[0,b_-]$, with $\Gamma(0)=0$. 
we define $b^*_e(u^0)>0$ (the threshold appearing in Theorem~\ref{th:cv exp}~\rm{(ii)}) as the largest value in $[0,b_-]$ such that
\begin{equation}\label{eq:def gamma}
\frac{\bigl(N_\infty^b\,b\bigr)^{2}\bigl(1+\mathcal{I}_b(0)\bigr)}{\varepsilon_b} <1\quad \text{ for all } b\in[0,b^*_e(u^0)).
\end{equation}

\begin{proof}[Proof of Theorem~\ref{th:cv exp}~\rm(ii)]
Fix $b < b^*_e(u^0)$, so that 
$(N^b_\infty b)^2\bigl(1+\mathcal{I}_b(0)\bigr)<\varepsilon_b$, and let
\[
\varepsilon:=\Bigl(\varepsilon_b\,(N_\infty^b b)^2\bigl(1+\mathcal{I}_b(0)\bigr)\Bigr)^{1/2}.
\]
By construction, we have 
\[
\begin{cases} \varepsilon < \varepsilon_b, \quad\text{hence}\quad \nu_b(\varepsilon)\geq \frac{1}{2\alpha_b}>0
\ \text{ by \eqref{eq:minoration nu}}, \\[6pt]
c_b:=\varepsilon - (N_\infty^b\,b)^2\bigl(1+\mathcal{I}_b(0)\bigr)>0 . \end{cases}
\]
(The case $b=0$ is linear: the bracket in \eqref{eq: I' joli} reduces to $\varepsilon$, and the argument
below applies with $\varepsilon=\varepsilon_0$.) Since $b<b_-$ and $\varepsilon\in(0,\varepsilon_b]$, inequality
\eqref{eq: I' joli} holds, and we show that its bracket
$\varepsilon-(N_\infty^b\,b)^2\bigl(1+\mathcal{I}_b(t)\bigr)$ stays positive for all $t\geq 0$.

From the entropy identity \eqref{eq:entropienonlin}, the a priori bound
\eqref{eq:a priori bound solution FPSL} and $N\in L^2_{\mathrm{loc}}(\R^+)$, we have, for any $t,t'\geq 0$,
\begin{equation}\label{eq: continuite I}
    \bigl|\mathcal{I}_b(t)-\mathcal{I}_b(t')\bigr|\leq \Big|\int_{t'}^{t}q(\tau)\,d\tau\Big|, \qquad q\in L^1_{\mathrm{loc}}(\R^+),
\end{equation}
so that $\mathcal{I}_b$ is continuous. Let
\[
t^* := \sup\bigl\{t\geq 0 \;:\; \mathcal{I}_b(s)\leq \mathcal{I}_b(0)\ \text{ for all } s\in[0,t]\bigr\},
\]
and assume by contradiction that $t^*<\infty$. By continuity $\mathcal{I}_b(t^*)\leq\mathcal{I}_b(0)$, hence
\[
(N_\infty^b\,b)^2\bigl(1+\mathcal{I}_b(t^*)\bigr)\leq (N_\infty^b\,b)^2\bigl(1+\mathcal{I}_b(0)\bigr)< \varepsilon,
\]
and, by \eqref{eq: continuite I}, there exists $\delta>0$ such that
$(N_\infty^b\,b)^2\bigl(1+\mathcal{I}_b(t)\bigr)<\varepsilon$ for all $t\in[t^*,t^*+\delta]$. On this interval the
bracket in \eqref{eq: I' joli} is non-negative, so $\mathcal{I}_b'(t)\leq -\nu_b(\varepsilon)\,\mathcal{I}_b(t)\leq 0$, and
therefore $\mathcal{I}_b(t)\leq \mathcal{I}_b(t^*)\leq \mathcal{I}_b(0)$ for all $t\in[t^*,t^*+\delta]$. This
contradicts the definition of $t^*$, hence $t^*=\infty$ and $\mathcal{I}_b(t)\leq\mathcal{I}_b(0)$ for all $t\geq0$.

Consequently the bracket in \eqref{eq: I' joli} is bounded below by $c_b$ for every $t\geq0$, and, setting
$\lambda_b:=\nu_b(\varepsilon)$,
\begin{equation}\label{eq: I excit final}
\fa t\geq 0,\qquad
\mathcal{I}_b'(t)+\lambda_b\,\mathcal{I}_b(t)\leq - c_b\Big(\frac{N(t)}{N_\infty^b}-1\Big)^{\!2}\leq 0 .
\end{equation}
By Grönwall's lemma, we conclude that for any $t \geq 0$,
\[
\mathcal{I}_b(t)
= \int_\mathbb{R}\!\left(\frac{u(t)}{u_\infty^b}-1\right)^{\!2}u_\infty^b\,dx
\leq e^{-\lambda_b\,t}
\int_\mathbb{R}\!\left(\frac{u^0}{u_\infty^b}-1\right)^{\!2}u_\infty^b\,dx,
\]
which proves the first inequality in \eqref{eq:cv exp u excit}.
Finally, multiplying \eqref{eq: I excit final} by $e^{\lambda_bt}$ and integrating over $(0,\infty)$, we get
\[
c_b\int_0^{\infty}\!\left(\frac{N(t)}{N_\infty^b}-1\right)^{\!2}e^{\lambda_bt}\,dt
\leq \mathcal{I}_b(0)<\infty,
\]
which proves the second inequality in \eqref{eq:cv exp u excit}.
\end{proof}

\subsection{$L \log L$-Cesàro bound for excitatory connection}

The major difficulty comes from the slow decay of $u_\infty$ at $+\infty$ which prevents the use of log-Sobolev inequalities. We are obliged to devise an adapted analysis.

\begin{prop}\label{prop:ces LlogL}
    For any $ b\geq 0$ consider $u$ solution of \eqref{eq:FPSL Non-lin} with $u^0\in L^2_w(\R)$ and $w(t)=\frac{u(t)}{u_\infty}$, then, we have
    \begin{equation}
        \sup_{T\geq 0}\frac{1}{T+1}\Bigg(\int_\R w(T)\log\big(w(T)\big)u_\infty\,dx +\int_0^T\int_\R\big(\partial_x\sqrt{w(t)}\big)^2u_\infty\,dx\,dt+\int_0^T\frac{N(t)}{N_\infty}\log\frac{N(t)}{N_\infty}\,dt\Bigg) <\infty.
    \end{equation}
    
\end{prop}

\begin{proof}
We consider the entropy equation \eqref{eq:entropienonlin} with the convex function $H(x) = x\log(x)$. Since $xH'(x)-H(x)=x$, this yields
\begin{equation*}
     \frac{d}{dt}\int_\mathbb{R} w(t)\log w(t)\, u_\infty\,dx = -\mathrm{Jump}(t)
     -4\int_\mathbb{R}\bigl(\partial_x\sqrt{w(t)}\bigr)^2 u_\infty\,dx
           - bN(t)\int_\mathbb{R} w(t)\,u_\infty'\,dx.
\end{equation*}
where the entropy dissipation due to jumps is defined as 
\begin{align*}
    \mathrm{Jump}(t) &= \frac{N(t)}{N_\infty}\log\frac{N(t)}{N_\infty}
           - w(t,0)\log w(t,0)
           - \bigl(\log w(t,0)+1\bigr)\Bigl(\frac{N(t)}{N_\infty}-w(t,0)\Bigr)\\
           &=\frac{N(t)}{N_\infty}\log\frac{N(t)}{N_\infty} + w(t,0)-\frac{N(t)}{N_\infty}\big(1+\log w(t,0)\big) \geq 0.
\end{align*}
Using that $-u_\infty'\leq \kappa_\infty u_\infty$ and $b>0$, we get
\[
-bN(t)\int_\R w(t)\,u_\infty'\,dx\leq \kappa_\infty bN(t)\int_\R w(t) u_\infty\,dx = \kappa_\infty bN(t).
\]

To handle the jump term of the entropy equality, we use the Young-Fenchel inequality \cite[p.11]{brezisFunctionalAnalysisSobolev2011} of the conjugated-convex couple $(x\log x, e^{y-1})$. For any $x,y>0$ and $0<\varepsilon <1$, we have
\[
xy\leq \varepsilon x\log(\varepsilon x) +e^{\frac{y}{\varepsilon}-1}.
\]
Therefore, taking $x= \frac{N(t)}{N_\infty}$ and $y=1+\log w(t,0)$, we obtain
\begin{equation}\label{eq:minoration jump}
    \mathrm{Jump}(t)\geq (1-\varepsilon)\frac{N(t)}{N_\infty}\log\frac{N(t)}{N_\infty}-e^{\frac{1}{\varepsilon}-1}w(t,0)^{\frac{1}{\varepsilon}}.
\end{equation}
In particular, for $\varepsilon=\frac{2}{3}$, we establish below that 
\begin{equation}\label{eq:bborne L^inf H^1}
     \exists K>0,\quad \fa t\geq 0, \quad \sqrt{e}w(t,0)^{\frac{3}{2}}\leq 3 \int_\R \big(\partial_x\sqrt{w(t,x)}\big)^2\,u_\infty\,dx +K.
\end{equation}
Adding \eqref{eq:minoration jump} and \eqref{eq:bborne L^inf H^1} with the entropy equality ensures that 
\[
\frac{d}{dt}\int_\R w(t)\log w(t)\, u_\infty\,dx + \int_\R \big(\partial_x\sqrt{w(t,x)}\big)^2\,u_\infty\,dx+\frac{1}{3}\frac{N(t)}{N_\infty}\log\frac{N(t)}{N_\infty}\leq K+\kappa_\infty bN(t).
\]
Integrating over $(0,T)$ and using \eqref{eq:est L^1 N}, we obtain,
\begin{equation}\label{eq:LlogL integree}
\int_\mathbb{R} w(T)\log w(T)\, u_\infty\,dx
+\int_0^T\!\!\int_\R\big(\partial_x\sqrt{w(t)}\big)^2u_\infty\,dx\,dt
+\frac{1}{3}\int_0^T\frac{N(t)}{N_\infty}\log\frac{N(t)}{N_\infty}\,dt
\;\leq\;\mathcal H(u^0)+2b\kappa_\infty+\big(K+2b\kappa_\infty\gamma_b\big)T ,
\end{equation}
 where
\[
  \mathcal H(u^0):=\int_\R u^0\log \frac{u^0}{u_\infty}\,\,dx.
\]
Dividing
\eqref{eq:LlogL integree} by $T+1$ yields
\[
\sup_{T\geq 0}\frac{1}{T+1}\bigg[\int_\mathbb{R} w(T)\log w(T)\, u_\infty\,dx+\int_0^T\int_\R\big(\partial_x\sqrt{w(t)}\big)^2u_\infty\,dx\,dt+\frac{1}{3}\int_0^T\frac{N(t)}{N_\infty}\log\frac{N(t)}{N_\infty}\,dt \bigg]<\infty.
\]

It remains to establish \eqref{eq:bborne L^inf H^1}.
We set $m:=\big\|u_\infty^{-1}\big\|_{\infty,(0,1)}$ and
$D(t):=\int_\R\big(\partial_x\sqrt{w(t)}\big)^2u_\infty\,dx$. Let $v:=\sqrt{w(t)}$. For
$y\in(0,1)$, $v(0)^2=v(y)^2-2\int_0^y vv'$, so that, averaging in $y$ over $(0,1)$ and using the
Cauchy--Schwarz inequality,
\[
w(t,0)=v(0)^2\leq \int_{0}^1v^2\,dx+2\Big(\int_{0}^1v^2\,dx\Big)^{1/2}\Big(\int_{0}^1v'^2\,dx\Big)^{1/2}
\leq  m+2m\sqrt{D(t)} ,
\]
since $\int_{0}^1v^2\,dx=\int_{0}^1\frac{u}{u_\infty}\,dx\leq m$ and
$\int_{0}^1v'^2\,dx\leq m\,D(t)$. By convexity, $(a+b)^{3/2}\leq\sqrt2\,(a^{3/2}+b^{3/2})$, hence
\[
\sqrt{e}\,w(t,0)^{3/2}\leq\sqrt{2e}\,m^{3/2}+4\sqrt e\,m^{3/2}D(t)^{3/4}\leq 3D(t)+K,
\qquad K:=\sqrt{2e}\,m^{3/2}+e^2m^6,
\]
where we used $c\,s^{3/4}\leq 3s+\frac{c^4}{256}$ for $c,s\geq0$, with $c=4\sqrt e\,m^{3/2}$.
This shows \eqref{eq:bborne L^inf H^1} and concludes the proof.

\end{proof}

\begin{remark}
    Theorem \ref{th:cesaro bound}~(i) comes by identifying the terms in \eqref{eq: ces L^2 inhib}, which are
    \[
  C_1(u^0)=N_\infty\int_{\R}\Big(\frac{u^0}{u_\infty}\Big)^{2}u_\infty\,\varphi\,dx,
  \qquad C_2= N_\infty C_\infty ,
    \]
    and Theorem \ref{th:cesaro bound}~(ii) is obtained with \eqref{eq:LlogL integree}, that is 
    \[
  C_1'(b,u^0)=3\bigl(\mathcal H(u^0)+2b \kappa_\infty\bigr),\qquad
  C_2'(b)=3\bigl(K+2b \kappa_\infty \gamma_b\bigr).
    \]
\end{remark}

\subsection{Well-posedness of \eqref{eq:FPSL Non-lin} for an initial datum in $L^2_{\mathrm w}(\R)$}
\label{sec:existence L2w suite}

As announced in Section~\ref{sec:existenceetunicite}, we can now complete the proof of
Theorem~\ref{th:existence and uniqueness FPSLNL} by removing the restriction $u^0\in\mathcal{X}$
made in Proposition~\ref{prop:existence unicite glob u^0 dans X}. The major step is to obtain a
bound on the activity which is uniform along an approximating sequence, that is, a bound
depending on the initial datum only through $\|u^0\|_{L^2_{\mathrm w}}$. It is provided by the estimates of the present section, and
the argument splits according to the sign of $b$: for $b\leq0$ the nonlinear term has a
favorable sign and the $L^2$-Cesàro bound applies directly, whereas for $b\geq0$ one closes a
Grönwall inequality using \eqref{eq: borne u'/u} together with the uniform $L^1$ estimate
\eqref{eq:est L^1 N}.

\begin{prop}[Uniform local bounds on the activity]\label{prop:borne L2 uniforme N}
Let $b\in\R$ and let $u^0\in\mathcal{X}$ be a probability density. Then the solution $(u,N)$ of
\eqref{eq:FPSL Non-lin} given by Proposition~\ref{prop:existence unicite glob u^0 dans X}
satisfies \eqref{eq:borne L2 uniforme N} for every $T>0$, with a function $\mathcal{K}$
depending only on $T$, $|b|$ and $\|u^0\|_{L^2_{\mathrm w}}$ and \eqref{eq:borne L2} . In particular the bound does not
involve $\|u^0/Q\|_\infty$.
\end{prop}

\begin{proof}

\noindent\emph{First case: $b\leq0$.} The $L^2$-Cesàro bound \eqref{eq: cesar L^2 inhib} of
Proposition~\ref{prop:cesaro inhib} gives directly
\[
  \int_0^TN(t)^2\,dt\;\leq\;N_\infty\Big(\|u^0\|_{L^2_{\mathrm w}}^{2}+C_\infty T\Big),
\]
since $\int_\R\big(\tfrac{u^0}{u_\infty}\big)^2u_\infty\varphi\,dx\leq\|u^0\|^2_{L^2_{\mathrm w}}$,
and the two remaining quantities in \eqref{eq:borne L2 uniforme N} are then controlled by the a
priori bound \eqref{eq:a priori bound solution FPSL} taken with $J(t)=bN(t)$. All the constants involved depend on $u^0$
only through $\|u^0\|_{L^2_{\mathrm w}}$.

\medskip
\noindent\emph{Second case: $b\geq0$.} 
We use the entropy identity \eqref{eq:entropienonlin} with $H(x)=x^2$ and $b'=0$. Since
$H''\equiv2$ and $xH'(x)-H(x)=x^2$, and since the jump bracket reduces to a square, it reads,
with $I(t)=\int_\R w(t)^2u_\infty\,dx$,
\begin{equation}\label{eq:identite I nonlin}
  I'(t)+2\int_\R\big(\partial_xw(t)\big)^2u_\infty\,dx
  +N_\infty\Big(\frac{N(t)}{N_\infty}-w(t,0)\Big)^{2}
  \;=\;b\,N(t)\int_\R w(t)^2\big(-u_\infty'\big)\,dx .
\end{equation}

Here $bN(t)\geq0$, so that \eqref{eq: borne u'/u} gives
$b\,N(t)\int_\R w^2(-u_\infty')\,dx\leq\kappa_\infty\,b\,N(t)\,I(t)$ and, discarding the two
non-negative terms of the left-hand side of \eqref{eq:identite I nonlin},
\[
  I'(t)\;\leq\;\kappa_\infty\,b\,N(t)\,I(t).
\]
By Grönwall's lemma and the uniform $L^1$ estimate \eqref{eq:est L^1 N}, namely
$\int_0^{T}N(t)\,dt\leq2(\gamma_bT+1)$, we obtain
\begin{equation}\label{eq:borne sup I nonlin}
  \sup_{t\in[0,T]}I(t)\;\leq\;\mathcal{M}(T):=\|u^0\|_{L^2_{\mathrm w}}^{2}\,
  e^{2\kappa_\infty b(\gamma_bT+1)} ,
\end{equation}
which is the first term of \eqref{eq:borne L2 uniforme N}. Integrating now
\eqref{eq:identite I nonlin} over $(0,T)$ and using \eqref{eq:borne sup I nonlin} on the
right-hand side,
\begin{equation}\label{eq:borne dissipation nonlin}
  2\int_0^T\!\!\int_\R\big(\partial_xw\big)^2u_\infty\,dx\,dt
  +N_\infty\int_0^T\Big(\frac{N(t)}{N_\infty}-w(t,0)\Big)^{2}dt
  \;\leq\;\|u^0\|_{L^2_{\mathrm w}}^{2}+2\kappa_\infty b\,\mathcal{M}(T)\big(\gamma_bT+1\big).
\end{equation}
It remains to convert the second term of \eqref{eq:borne dissipation nonlin} into a bound on $N$
itself. Writing, for $y\in(0,1)$, $w(t,0)=w(t,y)-\int_0^y\partial_xw(t,z)\,dz$ and integrating in
$y$, we get, with $c_1=\big\|\tfrac{1}{u_\infty}\big\|_{\infty,(0,1)}$,
\[
  |w(t,0)|\;\leq\;c_1+\Big(c_1\int_0^1\big(\partial_xw(t,z)\big)^2u_\infty\,dz\Big)^{1/2},
\]
so that $\int_0^Tw(t,0)^2dt$ is bounded by \eqref{eq:borne dissipation nonlin}. Since
\[
  \Big(\frac{N(t)}{N_\infty}\Big)^{2}\leq2\Big(\frac{N(t)}{N_\infty}-w(t,0)\Big)^{2}+2\,w(t,0)^2 ,
\]
the bound \eqref{eq:borne L2 uniforme N} follows, again with constants depending on $u^0$ only
through $\|u^0\|_{L^2_{\mathrm w}}$.

To show \eqref{eq:borne L2}, we use that for any function $f$, for a restrictive case of the Sobolev injections (see \cite{brezisFunctionalAnalysisSobolev2011}), we have
\[
\|f\|_2^2\leq \|f\|_\infty \|f\|_1 \quad \text{and}\quad \|f\|_\infty^2\leq 2\|f'\|_2\|f\|_2.
\]
Therefore, we obtain
\[
\|f\|_\infty^3\leq 4\|f'\|_2^2\|f\|_1.
\]
By choosing $f=u(t)$, in order to obtain \eqref{eq:borne L2}, it suffices to show that $\partial_xu(t) \in L^2_{loc}(\R^+;L^2(\R))$.
We have
\begin{align*}
    \|\partial_x u(t)\|^2_2 
    &=\int_\R \big(\partial_x w(t) u_\infty +u_\infty' w(t)\big)^2\,dx\leq 2\int_\R \big(\partial_x w(t)\big)^2u_\infty^2 +\big(u_\infty')^2w^2\,dx\\
    &\leq 2\|u_\infty\|_\infty \underbrace{\int_\R \big(\partial_x w(t)\big)^2u_\infty\,dx}_{\in L^1_{loc}(\R^+) \:\:\text{by}\:\eqref{eq:borne dissipation nonlin}} +2\Big\|\frac{(u_\infty')^2}{u_\infty}\Big\|_\infty\underbrace{\int_\R w(t)^2u_\infty\,dx}_{\in L^\infty_{loc}(\R^+)}.
\end{align*}
Hence, \eqref{eq:borne L2} is proven.
\end{proof}

We are now in a position to complete the proof of the main theorem.

\begin{proof}[Proof of Theorem~\ref{th:existence and uniqueness FPSLNL}]
Uniqueness for an arbitrary datum $u^0\in L^2_{\mathrm w}(\R)$ is
Proposition~\ref{prop:unicite L2w}, and existence has been proved in
Proposition~\ref{prop:existence unicite glob u^0 dans X} when $u^0\in\mathcal{X}$. There only
remains to construct a solution for an arbitrary probability density
$u^0\in L^2_{\mathrm w}(\R)$.

\medskip
Clearly, we have $\mathcal{C}^\infty_c(\R)\subset \mathcal{X}$ and $\mathcal{C}^\infty_c(\R)$ is dense in $\big(L^2_w(\R),\|.\|_{L^2_w}\big)$. Let therefore $(u^0_n)_n \in \mathcal{X}^{\mathbb{N}}$ be a sequence of probability
densities with $\|u^0_n-u^0\|_{L^2_{\mathrm w}}\to0$, and let $(u_n,N_n)$ be the corresponding
global solutions given by Proposition~\ref{prop:existence unicite glob u^0 dans X}. Fix~$T>0$.
Since $\|u^0_n\|_{L^2_{\mathrm w}}$ is bounded, Proposition~\ref{prop:borne L2 uniforme N}
provides a constant $\mathcal K_T$, independent of~$n$, such that
\begin{equation}\label{eq:bornes uniformes n}
  \sup_n\Big[\sup_{[0,T]}\|u_n(t)\|^2_{L^2_{\mathrm w}}
  +\int_0^T\!\!\int_\R\big(\partial_xw_n\big)^2u_\infty\,dx\,dt+\int_0^TN_n^2\,dt\Big]
  \;\leq\;\mathcal K_T .
\end{equation}

\medskip
Fix $M>0$. On $(-M,M)$ the weight $u_\infty$ is bounded from above and from below by positive
constants, so \eqref{eq:bornes uniformes n} makes $(u_n)_n$ bounded in
$L^2\bigl(0,T;H^1(-M,M)\bigr)$. Moreover, reading the equation
\[
  \partial_tu_n=-\partial_x\big((h+bN_n)u_n\big)+\partial_x^2u_n+N_n\delta_0
\]
on $(-M,M)$, where $h$ is bounded, each term on the right-hand side is bounded in
$L^2\bigl(0,T;H^{-1}(-M,M)\bigr)$ uniformly in $n$: the first two by
\eqref{eq:bornes uniformes n}, and the last one because $\delta_0\in H^{-s}(-M,M)$ for every
$s>\tfrac12$ and $(N_n)_n$ is bounded in $L^2(0,T)$. Since the embedding
$H^1(-M,M)\hookrightarrow L^2(-M,M)$ is compact,
the Aubin--Lions--Simon lemma
\cite{simonCompactSetsSpace1987} yields, up to a subsequence, the strong convergence of $u_n$ in
$L^2\bigl(0,T;L^2(-M,M)\bigr)$. A diagonal extraction over $M\in\N$ gives
\[
  u_n\longrightarrow u \quad\text{strongly in } L^2\bigl(0,T;L^2_{\mathrm{loc}}(\R)\bigr),
  \qquad
  N_n\rightharpoonup N \quad\text{weakly in } L^2(0,T),
\]
and, by \eqref{eq:bornes uniformes n} and weak-$*$ compactness, $u$ belongs to
$L^\infty\bigl(0,T;L^2_{\mathrm w}(\R)\bigr)\cap L^2\bigl(0,T;H^1_{\mathrm w}(\R)\bigr)$.

\medskip
Let $\phi\in\mathcal{C}^\infty_c\bigl([0,T)\times\R\bigr)$. All the terms of the weak
formulation \eqref{eq:weak formulation FPS} are linear in $u_n$ or in $N_n$, except
$\iint bN_nu_n\,\partial_x\phi$. Since $\phi$ has compact support in $x$, this term is the
product of $N_n$, which converges weakly in $L^2(0,T)$, by $u_n\partial_x\phi$, which converges
strongly in $L^2\bigl(0,T;L^2(-M,M)\bigr)$ for every $M$ with
$\operatorname{supp}\phi\subset[0,T)\times(-M,M)$. It therefore converges to $\iint bNu\,\partial_x\phi$.
Passing to the limit in \eqref{eq:weak formulation FPS} shows that $(u,N)$ satisfies the weak
formulation with initial datum $u^0$.

\medskip

It remains to check that $N$ is indeed the outgoing flux of $u$ in the weak sense which, as
recalled in \eqref{eq:conservation de la masse}, amounts to the conservation of mass. This
follows from the tightness of the family $(u_n)_n$: by
Cauchy--Schwarz and \eqref{eq:bornes uniformes n},
\[
  \int_{|x|>M}u_n(t,x)\,dx\;\leq\;
  \Big(\int_\R\frac{u_n(t,x)^2}{u_\infty(x)}\,dx\Big)^{1/2}
  \Big(\int_{|x|>M}u_\infty(x)\,dx\Big)^{1/2}
  \;\leq\;\mathcal K_T^{1/2}\Big(\int_{|x|>M}u_\infty\Big)^{1/2},
\]
which tends to $0$ as $M\to\infty$, uniformly in $n$ and in $t\in[0,T]$, since $u_\infty$ is a
probability density. Hence $\int_\R u(t)\,dx=1$ for every $t$, and
$N(t)=\lim_{x\to+\infty}h(x)u(t,x)$.

\medskip
Therefore $(u,N)$ is a solution of \eqref{eq:FPSL Non-lin} on $(0,T)$ with datum $u^0$, and
Proposition~\ref{prop:unicite L2w} guarantees that it does not depend on the extracted
subsequence. As $T>0$ was arbitrary, the solution is global, and the a priori bound
\eqref{eq:a priori bound solution FPSL} passes to the limit by lower semicontinuity of the norms.
Notice that \eqref{eq:borne L2} also passes to the limit, so it is valid for every $u^0\in L^2_w(\R)$.
\end{proof}

\section{Conclusion}

Combining new Cesàro bounds for the activity with precise estimates for the linear problem with a
time-dependent drift, we prove global existence for every connectivity. This is in
sharp contrast with the finite-time blow-up that can occur when a firing threshold is imposed.
Moreover, for initial data in $\mathcal{X}$, the activity remains locally bounded  in time. In the
weakly connected regime, the solution converges exponentially towards its unique stationary state.

Several limitations of the present analysis suggest natural directions for future work. First, our
assumptions on the drift $h$ cover a broad class of polynomial and exponential growth rates, but
still exclude faster-growing drifts, such as $h(x)=e^{x^2}$ as $x\to+\infty$. In this setting, it is
not known whether solutions remain global for every $b\in\mathbb{R}$, even though stationary states
still exist.

A second open question concerns the relation between the present limiting model and the NNLIF model
with a finite firing threshold $V_F$. From a physical perspective, it would be important to prove
that, for the NNLIF model with drift $h+bN$, the limit $V_F\to+\infty$ is precisely described by the
solution constructed here. This would establish that the limiting dynamics do not depend on the
particular choice of a sufficiently large firing threshold, and could provide another viewpoint on
the disappearance of the blow-up mechanism in the limiting problem.

Another natural extension concerns the way network activity feeds back on each neuron. In
\eqref{eq:FPSL Non-lin}, the mean synaptic current is $bN(t)$, which assumes an activity-independent
synaptic efficacy. Short-term synaptic plasticity makes this efficacy depend on the firing rate: in
the mean-field limit of the Tsodyks--Markram model~\cite{tsodyksNeuralNetworksDynamic1998}, the
linear case corresponds to frequency-independent transmission, whereas facilitation makes the
synaptic current grow like $N^2$ at low firing rates and depression makes it saturate at high rates.
This suggests replacing $bN(t)$ by $b\,f\big(N(t)\big)$ for a general function $f$. In the
excitatory regime, the global existence proved here should not persist for sufficiently superlinear
couplings: for $f(N)=N^2$ and exponential drifts, we expect finite-time blow-up as soon as $b$ is
large enough, whatever the initial datum. This question will be addressed in a forthcoming work.

Finally, the long-time behavior raises several questions. In the excitatory regime, numerical
simulations consistently display convergence towards a stationary state for drifts of the form $h(x) = -x+\beta+e^{x-x_0}$, even when several
stationary states coexist, see Section~\ref{sec:stationary states}. However, it is possible to construct specific functions $h$ such that solutions converge to periodic functions in the excitatory regime; this will be addressed in a forthcoming work. 

In the strongly inhibitory
regime, the situation is interesting as well: numerical simulations reveal the
emergence of periodic solutions even without delay. This is in contrast with the NNLIF model, for which, to our knowledge, periodic
solutions are only observed when a transmission delay is introduced. A possible explanation is that
the time needed to reach $+\infty$ acts as an intrinsic delay. A theoretical study of these periodic solutions, and in particular
of the role played by the growth of $h$ at $+\infty$, is a natural continuation of this work. As for
the NNLIF model, introducing a transmission delay through a drift of the form $h+bN(t-d)$ also leads
to periodic solutions when $b$ is sufficiently negative. Understanding the existence and stability
of such periodic solutions, with or without delay, remains largely open and connects with the
questions investigated in~\cite{brunelDynamicsSparselyConnected, brunelFastGlobalOscillations1999,
ikedaTheoreticalStudyEmergence2022,
caceresGlobalintimeClassicalSolutions2018,
rieutordAsymptoticDynamicsInhibitory2026, caceresSequencePseudoequilibriaDescribes2024}.\\

\textbf{Acknowledgment.} This work was supported by the Fondation Simone et Cino Del Duca, Institut de France. 

%---------------------------------
\appendix
\section{SLIF model with a time-dependent drift $J(t)$}
\label{appendix}

We prove Theorem~\ref{th:existence unicite avec J(t)}. The comparison bounds
\eqref{est:borne L^infini J} and \eqref{est:borne L^1 J} have already been proved in
Section~\ref{sec:linear}. We first establish here the two remaining entropy estimates
announced there, and we then construct the solution using a piecewise constant input $J(t)$.

\begin{prop}[A priori estimates on $(u,N)$]\label{prop:estimation a priori u_J}
Let $(u,N)$ be the solution of \eqref{eq:FPSL with J(t)} with $J\in L^2_{loc}(\R^+)$. The bounds
\eqref{eq:a priori bound solution FPSL} and \eqref{est: borne sur N}, announced in
Theorem~\ref{th:existence unicite avec J(t)}, are satisfied.
\end{prop}

\begin{proof}
Write, for all $t\geq0$,
\[
 w(t)=\frac{u(t)}{u_\infty} \qquad\text{and}\qquad I(t) = \int_\R w(t,x)^2u_\infty(x)\,dx .
\]
We use the entropy equality \eqref{RE:generique} with the convex function $H(w)=w^2$ and
$R=- J(t)\,\partial_x u= -J(t)\partial_x\big(w(t)u_\infty\big)$ in \eqref{eq:generique}, which gives, for any $t\geq 0$,
\[
    \frac{d}{dt}I(t)\leq -2\int_\R\big(\partial_x w(t)\big)^2 u_\infty\,dx
   - N_\infty \Big(\frac{N(t)}{N_\infty}-w(t,0) \Big)^2 -2J(t)\int_\R w(t)\,\partial_x\big(w(t)u_\infty\big)\,dx .
\]
Bounding the last term by $2ab\leq a^2+b^2$, with $a=\big(\partial_xw\big)\sqrt{u_\infty}$ and
$b=|J|\,w\sqrt{u_\infty}$, absorbs half of the dissipation and leaves the single differential
inequality
\begin{equation}\label{est:entropie L2}
    \frac{d}{dt}I(t)+\int_\R\big(\partial_x w(t)\big)^2 u_\infty\,dx
    + N_\infty \Big(\frac{N(t)}{N_\infty}-w(t,0) \Big)^2 \;\leq\; J(t)^2 I(t) ,
\end{equation}
from which both estimates follow. Indeed, discarding the two non-negative terms of the left-hand
side and applying Gr\"onwall's lemma gives $I(t)\leq I(0)e^{\mathcal{J}(t)}$; integrating
\eqref{est:entropie L2} over $(0,T)$ and inserting this bound then yields, since
$\int_0^TJ(t)^2e^{\mathcal{J}(t)}\,dt=e^{\mathcal{J}(T)}-1$,
\begin{equation}\label{est:entropie L2 integree}
    I(T)+\int_0^T\!\!\int_\R\big(\partial_x w\big)^2 u_\infty\,dx\,dt
    + N_\infty \int_0^T\Big(\frac{N(t)}{N_\infty}-w(t,0) \Big)^2 dt
    \;\leq\; I(0)+\int_0^TJ(t)^2I(t)\,dt\;\leq\; I(0)\,e^{\mathcal{J}(T)} .
\end{equation}
Its first two terms are exactly \eqref{eq:a priori bound solution FPSL}.

It remains to convert the third one into a bound on $N$. For all $y\in (0,1)$ we have
\[
w(t,0)= w(t,y) - \int_0^y \partial_x w(t,z)\,dz
= \int_0^1 \frac{u(t,y)}{u_\infty (y)}\, dy -\int_0^1\!\!\int_0^y \partial_x w(t,z)\,dz\, dy ,
\]
which implies, with $c_1= \frac{1}{\min_{(0,1)} u_\infty }$ and using $|\int_\R u(t)\,dx|\leq \|u^0\|_1$,
\[
| w(t,0) | \leq c_1\|u^0\|_1 + \left(c_1 \int_0^1 (\partial_x w(t,z))^2 u_\infty\, dz \right)^{1/2}
\qquad\text{hence}\qquad
w(t,0)^2\leq 2(c_1\|u^0\|_1)^2+2c_1\int_\R \big(\partial_x w(t)\big)^2 u_\infty\, dz .
\]
Since $\big(\frac{N(t)}{N_\infty}\big)^2\leq 2\big(\frac{N(t)}{N_\infty}-w(t,0)\big)^2+2\,w(t,0)^2$,
integrating in time and using \eqref{est:entropie L2 integree} for both terms gives
\eqref{est: borne sur N}.
\end{proof}

Having the a priori estimates at hand, we can now prove the existence in
Theorem~\ref{th:existence unicite avec J(t)}.

\begin{proof}[Proof of Theorem \ref{th:existence unicite avec J(t)}]
First, for $J\in\R$, let $\big(S(J,t)\big)_{t\geq 0}$ denote the continuous semigroup associated with the Fokker--Planck equation
\begin{equation}\label{eq:FPSL avec J}
    \left\{
    \begin{aligned}
        &\frac{\partial u}{\partial t}(t,x)+\frac{\partial}{\partial x}\big(\big(h(x)+J\big)u(t,x)\big)-\frac{\partial^2 u}{\partial x^2}(t,x)=\delta_0(x)\,N(t),\quad x\in\R,\ t\geq 0,\\
        &N(t)=\lim_{x\to+\infty}h(x)u(t,x).
    \end{aligned}\right.
\end{equation}
Now let $T>0$. To prove existence, we split $[0,T]$ into small intervals on which the drift $J(t)$ is constant; the a priori estimates then allow us to pass to the limit. For a given $n\in\N$, we define the following objects:
\begin{enumerate}
    \item For all $k\in\{0,\dots,n\}$, $\displaystyle T_k^n=\frac{k}{n}T$.
    \item $\displaystyle J^{(n)}(t)=\sum_{k=0}^{n-1}J_k^n\,\mathbf{1}_{(T_k^n,T_{k+1}^n)}(t)$,\quad with $\displaystyle J_k^n=\frac{n}{T}\int_{T_k^n}^{T_{k+1}^n}J(\tau)\,d\tau$.
    \item For all $k$ and $t\in(T_k^n,T_{k+1}^n)$,
    \[
        u^{(n)}(t)=S(J_k^n,t-T_k^n)\,u^{(n)}(T_k^n),\quad
        w^{(n)}(t)=\frac{u^{(n)}(t)}{u_\infty},\quad
        N^{(n)}(t)=\lim_{x\to\infty}h(x)u^{(n)}(t,x).
    \]
    \item For all $n\in\N$ and $t\in(0,T)$,
    \[
        I^{(n)}(t)=\int_\R w^{(n)}(t)^2 u_\infty\,dx.
    \]
\end{enumerate}
The function $J^{(n)}$ is piecewise constant and, since $J\in L^2(0,T)$,
\[
    \lim_{n\to\infty}\int_0^T\big|J(t)-J^{(n)}(t)\big|^2\,dt=0.
\]
Applying \eqref{est:entropie L2 integree} on each interval $(T_k^n,T_{k+1}^n)$, where the input
is the constant $J_k^n$, and summing over $k$, the terms $I^{(n)}$ telescope. Chaining
Gr\"onwall's lemma in the same way gives $I^{(n)}(t)\leq I(0)\exp\big(\int_0^tJ^{(n)}(s)^2ds\big)$,
so that, exactly as above, for every $t\in(0,T)$,
\[
    I^{(n)}(t)+\int_0^t\!\!\int_\R\big(\partial_x w^{(n)}\big)^2u_\infty\,dx\,ds
    +N_\infty\int_0^t\Big(\frac{N^{(n)}(s)}{N_\infty}-w^{(n)}(s,0)\Big)^2ds
    \;\leq\;I(0)\exp\!\bigg(\int_0^TJ^{(n)}(s)^2\,ds\bigg)\;\leq\;I(0)\,e^{\mathcal{J}(T)},
\]
the last inequality following from Cauchy--Schwarz, since each $J_k^n$ is the mean of $J$ over
$(T_k^n,T_{k+1}^n)$. Arguing as in the proof of \eqref{est: borne sur N}, we conclude that
\begin{equation}\label{eq:borne unif I}
    \sup_{n\in\N}\bigg[\sup_{0\leq t\leq T}\int_\R w^{(n)}(t)^2u_\infty\,dx
    +\int_0^T N^{(n)}(t)^2\,dt\bigg]\;\leq\;C\big(T+e^{\mathcal{J}(T)}\big)\,I(0)\;<\;\infty .
\end{equation}

The sequence $\big(u^{(n)},N^{(n)}\big)_n$ is thus bounded in $L^\infty\big((0,T);L^2_{\mathrm w}(\R)\big)\times L^2(0,T)$. Hence, up to a subsequence, there exists $(u_T,N_T)\in L^\infty\big((0,T);L^2_{\mathrm w}(\R)\big)\times L^2(0,T)$ such that, in the weak-$*$ topology,
\[
    u^{(n)}\rightharpoonup u_T\quad\text{and}\quad N^{(n)}\rightharpoonup N_T.
\]

Let $\phi\in \mathcal{C}_c^\infty(\R^+\times\R)$ be a test function such that supp$(\phi)\subset [0,T)$. For every $n\in\N$, $k\in\{0,\dots,n\}$ and $u^{(n)}$ being a strong solution of \eqref{eq:FPSL avec J} with $J=J_k^n$ on $(T_k^n,T_{k+1}^n)$, it satisfies the weak formulation as follows
\begin{align*}
    \int_{T_k^n}^{T_{k+1}^n}\!\!\int_\R u^{(n)}(t,x)\Big(-\partial_t\phi(t,x)
    &-\big(h(x)+J_k^n\big)\partial_x\phi(t,x)-\partial_x^2\phi(t,x)\Big)\,dx\,dt\\
    &=\int_{T_k^n}^{T_{k+1}^n}N^{(n)}(t)\,\phi(t,0)\,dt
    +\int_\R u^{(n)}(T_k^n,x)\phi(T_k^n,x)\,dx\\
    &\quad -\int_\R u^{(n)}(T_{k+1}^n,x)\phi(T_{k+1}^n,x)\,dx.
\end{align*}
Summing over $k$, we obtain, for any test function $\phi$ with $\operatorname{supp}(\phi)\subset[0,T)$,
\begin{align*}
    \int_0^T\!\!\int_\R u^{(n)}(t,x)\Big(-\partial_t\phi(t,x)
    &-\big(h(x)+J^{(n)}(t)\big)\partial_x\phi(t,x)-\partial_x^2\phi(t,x)\Big)\,dx\,dt\\
    &=\int_0^T N^{(n)}(t)\,\phi(t,0)\,dt+\int_\R u^0(x)\phi(0,x)\,dx.
\end{align*}
Since $J^{(n)}\to J$ strongly, we may pass to the limit to obtain
\begin{align*}
    \int_0^T\!\!\int_\R u_T(t,x)\Big(-\partial_t\phi(t,x)
    &-\big(h(x)+J(t)\big)\partial_x\phi(t,x)-\partial_x^2\phi(t,x)\Big)\,dx\,dt\\
    &=\int_0^T N_T(t)\,\phi(t,0)\,dt+\int_\R u^0(x)\phi(0,x)\,dx.
\end{align*}

Therefore $(u_T,N_T)$ is a weak solution of \eqref{eq:FPSL avec J} on $(0,T)$.
The solution $u_T$ is unique. Indeed, if $(u^1_T,N^1_T)$ and $(u^2_T,N^2_T)$ are two weak
solutions of \eqref{eq:FPSL avec J} on $(0,T)$ with the same input $J$ and the same initial
datum, then, the problem being linear in the pair $(u,N)$, their difference is a weak solution
with zero initial datum. Applying \eqref{eq:a priori bound solution FPSL} to this difference gives
$\|u^1_T(t)-u^2_T(t)\|^2_{L^2_{\mathrm w}}\leq e^{\mathcal{J}(t)}\cdot0=0$ for every
$t\in(0,T)$, whence $u^1_T=u^2_T$; the activity, being the outgoing flux
$N_T(t)=\lim_{x\to+\infty}h(x)u_T(t,x)$, is then determined as well.
This allows to define the function globally on $\R^+$ by setting, for all $t\geq 0$,
\[
    u(t)=u_{\lfloor t\rfloor+1}(t),\qquad N(t)=N_{\lfloor t\rfloor+1}(t).
\]
By construction and by uniqueness of $u_T$, the pair $(u,N)$ is well defined and unique on $\R^+$, and it satisfies the weak formulation \eqref{eq:weak formulation FPS}. This concludes the proof of Theorem \ref{th:existence unicite avec J(t)}.
\end{proof}

\section*{Declarations}
\noindent\textbf{Conflict of interest.} The authors declare no conflict of interest.

\noindent\textbf{Data availability.} No data have been used or produced.

\noindent\textbf{AI-assisted tools.}
During the preparation of this work, the authors used a generative AI assistant
(\emph{Claude - Opus 5.5}, Anthropic) to improve the language of the manuscript, to check its bibliographic
data, and to carry out a systematic critical review of the statements, computations and internal
consistency of the proofs. Every suggestion was independently checked by the authors, who
conceived all the mathematical content and take full responsibility for the article.

\bibliographystyle{siam}  
\bibliography{References_Zotero_septembre_2026}

\begin{thebibliography}{10}

\bibitem{bretteAdaptiveExponentialIntegrateandfire2005}
{\sc R.~Brette and W.~Gerstner}, {\em Adaptive {Exponential}
  {Integrate}-and-{Fire} {Model} as an {Effective} {Description} of {Neuronal}
  {Activity}}, Journal of Neurophysiology, 94 (2005), pp.~3637--3642.

\bibitem{brezisFunctionalAnalysisSobolev2011}
{\sc H.~Brezis}, {\em Functional {Analysis}, {Sobolev} {Spaces} and {Partial}
  {Differential} {Equations}}, Springer New York, New York, NY, 2011.

\bibitem{brunelDynamicsSparselyConnected}
{\sc N.~Brunel}, {\em Dynamics of sparsely connected networks of excitatory and
  inhibitory spiking neurons}, Journal of Computational Neuroscience, 8 (2000),
  pp.~183--208.

\bibitem{brunelFastGlobalOscillations1999}
{\sc N.~Brunel and V.~Hakim}, {\em Fast global oscillations in networks of
  integrate-and-fire neurons with low firing rates}, Neural Computation, 11
  (1999), pp.~1621--1671.

\bibitem{brunelFiringRateNoisy2003}
{\sc N.~Brunel and P.~E. Latham}, {\em Firing {Rate} of the {Noisy} {Quadratic}
  {Integrate}-and-{Fire} {Neuron}}, Neural Computation, 15 (2003),
  pp.~2281--2306.

\bibitem{caceresSequencePseudoequilibriaDescribes2024}
{\sc M.~J. C{\'a}ceres, J.~A. Ca{\~n}izo, and A.~Ramos-Lora}, {\em Sequence of
  pseudoequilibria describes the long-time behavior of the nonlinear noisy
  leaky integrate-and-fire model with large delay}, Physical Review E, 110
  (2024), p.~064308.

\bibitem{caceresAnalysisNonlinearNoisy2010}
{\sc M.~J. C{\'a}ceres, J.~A. Carrillo, and B.~Perthame}, {\em Analysis of
  nonlinear noisy integrate \& fire neuron models: blow-up and steady states},
  The Journal of Mathematical Neuroscience, 1 (2011), p.~7.

\bibitem{caceresGlobalintimeClassicalSolutions2018}
{\sc M.~J. C{\'a}ceres, P.~Roux, D.~Salort, and R.~Schneider}, {\em
  Global-in-time classical solutions and qualitative properties for the {NNLIF}
  neuron model with synaptic delay}, Communications in Partial Differential
  Equations, 44 (2019), pp.~1358--1386.

\bibitem{carrilloQualitativePropertiesSolutions2015}
{\sc J.~Carrillo, B.~Perthame, D.~Salort, and D.~Smets}, {\em Qualitative
  {Properties} of {Solutions} for the {Noisy} {Integrate} \& {Fire} model in
  {Computational} {Neuroscience}}, Nonlinearity, 28 (2015), p.~3365.

\bibitem{carrilloClassicalSolutionsNonlinear2011}
{\sc J.~A. Carrillo, M.~d.~M. Gonz{\'a}lez, M.~P. Gualdani, and M.~E.
  Schonbek}, {\em Classical solutions for a nonlinear {Fokker}-{Planck}
  equation arising in computational neuroscience}, Communications in Partial
  Differential Equations, 38 (2013), pp.~385--409.

\bibitem{carrilloNonlinearPartialDifferential2025}
{\sc J.~A. Carrillo and P.~Roux}, {\em Nonlinear partial differential equations
  in neuroscience: from modelling to mathematical theory}, Mathematical Models
  and Methods in Applied Sciences, 35 (2025), pp.~403--584.

\bibitem{delarueGlobalSolvabilityNetworked2015}
{\sc F.~Delarue, J.~Inglis, S.~Rubenthaler, and E.~Tanr{\'e}}, {\em Global
  solvability of a networked integrate-and-fire model of {McKean}-{Vlasov}
  type}, The Annals of Applied Probability, 25 (2015).

\bibitem{delarueParticleSystemsSingular2015}
\leavevmode\vrule height 2pt depth -1.6pt width 23pt, {\em Particle systems
  with a singular mean-field self-excitation. {Application} to neuronal
  networks}, Stochastic Processes and their Applications, 125 (2015),
  pp.~2451--2492.

\bibitem{douNoisyIntegrateandfireEquation2024}
{\sc X.~Dou, B.~Perthame, D.~Salort, and Z.~Zhou}, {\em Noisy
  integrate-and-fire equation: continuation after blow-up}, Sept. 2024.
\newblock Preprint, arXiv:2409.14749.

\bibitem{zbMATH08050032}
{\sc X.~Dou and Z.~Zhou}, {\em Dilating blow-up time: a generalized solution of
  the {NNLIF} neuron model and its global well-posedness}, Commun. Math. Sci.,
  23 (2025), pp.~1081--1138.

\bibitem{fournierToyModelInteracting2015}
{\sc N.~Fournier and E.~L{\"o}cherbach}, {\em On a toy model of interacting
  neurons}, Annales de l'Institut Henri Poincar{\'e}, Probabilit{\'e}s et
  Statistiques, 52 (2016).

\bibitem{hamblyMcKeanVlasovEquationPositive2019}
{\sc B.~Hambly, S.~Ledger, and A.~S{\o}jmark}, {\em A {McKean}--{Vlasov}
  equation with positive feedback and blow-ups}, The Annals of Applied
  Probability, 29 (2019), pp.~2338--2373.

\bibitem{ikedaTheoreticalStudyEmergence2022}
{\sc K.~Ikeda, P.~Roux, D.~Salort, and D.~Smets}, {\em Theoretical study of the
  emergence of periodic solutions for the inhibitory {NNLIF} neuron model with
  synaptic delay}, Mathematical Neuroscience and Applications, Volume 2 (2022).

\bibitem{liuInvestigatingIntegrateFire2021}
{\sc J.-G. Liu, Z.~Wang, Y.~Xie, Y.~Zhang, and Z.~Zhou}, {\em Investigating the
  integrate and fire model as the limit of a random discharge model: a
  stochastic analysis perspective}, Mathematical Neuroscience and Applications,
  1 (2021), pp.~Art. 2, 1--36.

\bibitem{mattiaPopulationDynamicsInteracting2002}
{\sc M.~Mattia and P.~Del~Giudice}, {\em Population dynamics of interacting
  spiking neurons}, Physical Review E, 66 (2002), p.~051917.

\bibitem{michelGeneralRelativeEntropy2005}
{\sc P.~Michel, S.~Mischler, and B.~Perthame}, {\em General relative entropy
  inequality: an illustration on growth models}, Journal de Math{\'e}matiques
  Pures et Appliqu{\'e}es, 84 (2005), pp.~1235--1260.

\bibitem{montbrioMacroscopicDescriptionNetworks2015}
{\sc E.~Montbri{\'o}, D.~Paz{\'o}, and A.~Roxin}, {\em Macroscopic description
  for networks of spiking neurons}, Physical Review X, 5 (2015), p.~021028.

\bibitem{muckenhouptHardysInequalityWeights1972}
{\sc B.~Muckenhoupt}, {\em Hardy's inequality with weights}, Studia
  Mathematica, 44 (1972), pp.~31--38.

\bibitem{perthameTransportEquationsBiology2007}
{\sc B.~Perthame}, {\em Transport Equations in Biology}, Frontiers in
  Mathematics, Birkh{\"a}user, Basel, 2007.

\bibitem{perthameFokkerPlanckEquationSuperlinear2026}
{\sc B.~Perthame, C.~Rieutord, and D.~Salort}, {\em A {Fokker}-{Planck}
  equation with superlinear drift at infinity for {Integrate}-and-{Fire}
  model}, Jan. 2026.
\newblock Preprint, arXiv:2601.19282.

\bibitem{praMcKeanVlasovLimitInteracting1996}
{\sc P.~D. Pra and F.~D. Hollander}, {\em {McKean}-{Vlasov} limit for
  interacting random processes in random media}, Journal of Statistical
  Physics, 84 (1996), pp.~735--772.

\bibitem{ricciardiOrnsteinUhlenbeckProcessModel1979}
{\sc L.~M. Ricciardi and L.~Sacerdote}, {\em The {Ornstein}-{Uhlenbeck} process
  as a model for neuronal activity. {I}. {Mean} and variance of the firing
  time}, Biological Cybernetics, 35 (1979), pp.~1--9.

\bibitem{rieutordAsymptoticDynamicsInhibitory2026}
{\sc C.~Rieutord and D.~Salort}, {\em Asymptotic dynamics of inhibitory
  networks for the {NNLIF} {Model} in the large-delay limit}.
\newblock Preprint, June 2026.

\bibitem{simonCompactSetsSpace1987}
{\sc J.~Simon}, {\em Compact sets in the space {$L^p(0,T;B)$}}, Annali di
  Matematica Pura ed Applicata, 146 (1987), pp.~65--96.

\bibitem{tsodyksNeuralNetworksDynamic1998}
{\sc M.~Tsodyks, K.~Pawelzik, and H.~Markram}, {\em Neural networks with
  dynamic synapses}, Neural Computation, 10 (1998), pp.~821--835.

\end{thebibliography}

\end{document}